\documentclass{article}
\usepackage{mathrsfs}
\usepackage{amsmath}
\usepackage{amssymb, amsmath}
\usepackage{graphicx}
\usepackage{fancyhdr}
\usepackage{enumerate}
\usepackage{CJK}
\catcode`\@=11 \@addtoreset{equation}{section}

\catcode`\@=12

\usepackage{color}

\newtheorem{Theorem}{Theorem}[section]
\newtheorem{Proposition}{Proposition}[section]
\newtheorem{Lemma}{Lemma}[section]
\newtheorem{Corollary}{Corollary}[section]

\newtheorem{Remark}{Remark}[section]

\newcommand{\newcom}{\newcommand}
\newcommand{\bTheorem}[1]{
\begin{Theorem} \label{T#1} }
\newcommand{\eT}{\end{Theorem}}

\newcommand{\bProposition}[1]{
\begin{Proposition} \label{P#1}}
\newcommand{\eP}{\end{Proposition}}

\newcommand{\bLemma}[1]{
\begin{Lemma} \label{L#1} }
\newcommand{\eL}{\end{Lemma}}

\newcommand{\bCorollary}[1]{
\begin{Corollary} \label{C#1} }
\newcommand{\eC}{\end{Corollary}}

\newcommand{\beq}{\begin{equation}}
\newcommand{\eeq}{\end{equation}}
\newcom{\ben}{\begin{eqnarray}}
\newcom{\een}{\end{eqnarray}}
\newcom{\beno}{\begin{eqnarray*}}
\newcom{\eeno}{\end{eqnarray*}}
\newcom{\bali}{\begin{aligned}}
\newcom{\eali}{\end{aligned}}

\newcommand{\bFormula}[1]{
\begin{equation} \label{#1}}
\newcommand{\eF}{\end{equation}}

\newcommand{\f}{\frac}
\newcommand{\Om}{\Omega}

\newcommand{\p}{\partial}

\newcommand{\vr}{\varrho}

\newcommand{\vt}{\vartheta}
\newcommand{\vu}{\mathbf{u}}

\newcommand{\ma}{\mathcal{A}}

\newcommand{\vb}{\mathbf{B}}

\newcommand{\Div}{{\rm div}}
\newcommand{\Grad}{\nabla}

\newcommand{\dx}{{\rm d} x}

\newcommand{\dt}{{\rm d} t }
\newcommand{\ds}{{\rm d} s}

\newcommand{\ep}{\varepsilon}

\font\F=msbm10 scaled 1000
\newcommand{\R}{\mbox{\F R}}

\newcommand\Cbox[2]{%
    \newbox\contentbox%
    \newbox\bkgdbox%
    \setbox\contentbox\hbox to \hsize{%
        \vtop{
            \kern\columnsep
            \hbox to \hsize{%
                \kern\columnsep%
                \advance\hsize by -2\columnsep%
                \setlength{\textwidth}{\hsize}%
                \vbox{
                    \parskip=\baselineskip
                    \parindent=0bp
                    #2
                }%
                \kern\columnsep%
            }%
            \kern\columnsep%
        }%
    }%
    \setbox\bkgdbox\vbox{
        \color{#1}
        \hrule width  \wd\contentbox %
               height \ht\contentbox %
               depth  \dp\contentbox
        \color{black}
    }%
    \wd\bkgdbox=0bp%
    \vbox{\hbox to \hsize{\box\bkgdbox\box\contentbox}}%
    \vskip\baselineskip%
}

\begin{document}

%%%%%%%%%%%%%%%%%%%%%%%%%%%%%%%%

\pagestyle{fancy} \lhead{\color{blue}{Global well-posedness for 3D full non-resistive MHD}} \rhead{\emph{Y.Li}}

\title{\bf Global well-posedness for the three-dimensional full compressible viscous non-resistive MHD system }

\author{
Yang Li$^{\dag,\ddag}$
 \\ \\  $^{\dag}$School of Mathematical Sciences, \\
Anhui University, Hefei 230601, People's Republic of China \\ Email: lynjum@163.com \\ \\
$^\ddag$Center of Pure Mathematics, \\
Anhui University, Hefei 230601, People's Republic of China \\ \\
}

\maketitle
{ \centerline {\bf Abstract}}
In this paper, we consider the three-dimensional full compressible viscous non-resistive MHD system. Global well-posedness is proved for an initial-boundary value problem around a strong background magnetic field. It is also shown that the unique solution converges to the steady state at an almost exponential rate as time tends to infinity. The proof is based on the celebrated two-tier energy method, due to Guo and Tice [\emph{Arch. Ration. Mech. Anal.}, 207 (2013), pp. 459--531; \emph{Anal. PDE.}, 6 (2013), pp. 287--369], reformulated in Lagrangian coordinates. The obtained result may be viewed as an extension of Tan and Wang [\emph{SIAM J. Math. Anal.}, 50 (2018), pp. 1432--1470] to the context of heat-conductive fluids. This in particular verifies the stabilization effects of vertical magnetic field in the full compressible non-resistive fluids.
{
}



\tableofcontents

\section{Introduction}
\subsection{Background and motivation}\label{back}

The time evolution of electrically conductive fluids in the presence of magnetic field is described by the magnetohydrodynamic equations (`MHD' in short). It is widely applied in astrophysics, thermonuclear reactions and industry, among may others. Assuming that the fluids are compressible viscous and ideal polytropic, whereas neglecting the resistivity coefficient and heat conductivity coefficient, the governing equations read as the following system of partial differential equations in Eulerian coordinates (see \cite{HC,DF,HW2}):
\beq\label{fmhd_1}
\p_t \tilde{\vr} +\Div_y (\tilde{\vr} \tilde{\vu})=0,
\eeq
\beq\label{fmhd_2}
\tilde{\vr}(\p_t \tilde{\vu}+ \tilde{\vu}\cdot \Grad_y \tilde{\vu})-\mu \Delta_y \tilde{\vu}-(\mu+\lambda)\Grad_y \Div_y \tilde{\vu}
+\Grad_y \left(
\tilde{\vr} \tilde{\vt}+\f{1}{2}|\tilde{\vb}|^2
\right)=\tilde{\vb} \cdot \Grad_y \tilde{\vb},
\eeq
\beq\label{fmhd_3}
\p_t \tilde{\vb} + \tilde{\vu}\cdot \Grad_y \tilde{\vb} -\tilde{\vb} \cdot \Grad_y \tilde{\vu} +\tilde{\vb}\Div_y \tilde{\vu}=\mathbf{0},
\eeq
\beq\label{fmhd_4}
\tilde{\vr}(\p_t \tilde{\vt}+ \tilde{\vu}\cdot \Grad_y \tilde{\vt})+ \tilde{\vr} \tilde{\vt}  \Div_y \tilde{\vu}
=\kappa \Delta_y \tilde{\vt}+
2\mu |\mathbb{D}_y\tilde{\vu} |^2+\lambda (\Div_y \tilde{\vu} )^2
\eeq
\beq\label{fmhd_5}
\Div_y \tilde{\mathbf{B}}=0.
\eeq
Here the unknowns $\tilde{\vr}(t,y),\tilde{\vu}(t,y)\in\R^3,\tilde{\vb}(t,y)\in\R^3,\tilde{\vt}(t,y)$ represent the density, the velocity field, the magnetic field and the temperature respectively. $t\in \R^{+}$ and $y\in \Om$ are the time variable and space variable respectively. In this paper, we suppose $\Om$ to be an infinite slab $\Om:=\R^2\times (0,1)$. $\mu$ and $\lambda$ denote the shear and bulk viscosity coefficients respectively, obeying $\mu>0,\,2\mu+3\lambda\geq 0$. $\kappa>0$ is the heat-conductivity coefficient. $\mathbb{D}_y\tilde{\vu} $ signifies the symmetric gradient of velocity, meaning $\mathbb{D}_y\tilde{\vu}=\f{\Grad_y \tilde{\vu}+(\Grad_y \tilde{\vu})^t}{2} $.

The equations (\ref{fmhd_1})-(\ref{fmhd_5}) are supplemented with the initial conditions
\beq\label{fmhd_6}
(\tilde{\vr},\tilde{\vu},\tilde{\vb},\tilde{\vt})|_{t=0}=(\tilde{\vr}_0,\tilde{\vu}_0,\tilde{\vb}_0,\tilde{\vt}_0),
\eeq
as well as the boundary conditions
\begin{equation}\label{fmhd_7}
\left\{\begin{aligned}
& \tilde{\vu}|_{\p \Om}=\mathbf{0}, \\
&  \tilde{\vt}|_{\p \Om}=\overline{\vt},\\
\end{aligned}\right.
\end{equation}
for any given constant $\overline{\vt}>0$.

Below we review some previous results. There are many results on the incompressible MHD system, meaning $\Div_y \tilde{\vu}=0$. In analogy with the incompressible Navier-Stokes system, global well-posedness to the viscous resistive incompressible MHD system is classical. From the physical point of view, the Reynolds number and the electric conductivity may be large in several models. It is thus interesting and natural to consider the MHD system with only kinematic viscosity or with only finite electric conductivity or without any dissipation effects. In 3D, Bardos et al. \cite{BSS1} proved global well-posedness for inviscid non-resistive incompressible MHD system in the presence of a strong magnetic field; see also Cai and Lei \cite{CaLe1} as well as He et al. \cite{HXY1} for recent treatment. Cao and Wu \cite{CaWu1} proved the existence of weak solutions for 2D inviscid resistive incompressible MHD system. However, for viscous non-resistive fluids, the handling seems to be more tricky and thus received more attention in recent years. Lin et al. \cite{LXZ} proved global small solutions for 2D viscous non-resistive incompressible MHD system, see also \cite{RWZ,ZT} for more related results in 2D. In 3D, Xu and Zhang \cite{XZ} obtained the global small solutions for viscous non-resistive incompressible MHD system; relevant results were later established by Abidi and Zhang \cite{AbZh1}, Pan et al. \cite{PZZu1}, Tan and Wang \cite{TW} via different approaches. Very recently, Wu and Zhu \cite{WuZi1} proved global well-posedness for 3D incompressible MHD system with mixed partial dissipation and magnetic diffusion.

Let us come back to the compressible MHD system. Based on the compressible Navier-Stokes system \cite{LS,FNP}, Ducomet and Feireisl \cite{DF} proved the global existence of weak solutions to 3D full compressible viscous resistive MHD system, see also Hu and Wang \cite{HW2} for similar result. Kawashima \cite{Kaw1} proved global well-posedness for compressible viscous resistive MHD system around constant equilibrium states. Li et al. \cite{LiXuZw1} proved global classical solutions to 3D compressible viscous resistive MHD system with large oscillations and vacuum. Feireisl and Li \cite{FL1} proved the existence of infinitely many global weak solutions to 3D compressible barotropic inviscid resistive MHD system. It turns out that the problem becomes quite delicate for compressible viscous non-resistive MHD system. There are satisfactory results in the simplified 1D geometry. Jiang and Zhang \cite{JZW} proved the existence and uniqueness of global strong solution to the isentropic case with large initial data. We refer to \cite{LY2,LY3,LS1D} for more results in 1D concerning isentropic and heat-conductive non-resistive MHD system with large initial data. Wu and Wu \cite{WW} proved global well-posedness for 2D compressible isentropic viscous non-resistive MHD system by requiring smallness of initial data. Very recently, we considered a special\footnote{In our model we assumed that the motion of fluids takes place in the plane while the magnetic field acts on the fluids only in the vertical direction. This makes the non-resistive magnetic equation to be the continuity equation transported by the same velocity field like the density.} 2D compressible viscous non-resistive MHD system both for isentropic and heat-conductive fluids and obtained the existence of global weak solutions with large initial data \cite{LS2D,LS2Dfu}; see also \cite{LiZh1} for an extension including density-dependent viscosity coefficient and non-monotone pressure law. Tan and Wang \cite{TW} obtained global well-posedness for 3D compressible barotropic viscous non-resistive MHD system.

As far as the variation of temperature is taken into account, the only available results for compressible viscous non-resistive MHD system in multi-dimensions are \emph{locally} in time, see for instance \cite{FnYu1,Zh1}. Thus the \emph{global} well-posedness for \emph{full} compressible viscous non-resistive MHD system in multi-dimensions remains completely open up to now. The present paper aims to solve this problem partially. To do this, we prescribe the problem in an infinite slab $\R^2\times (0,1)$ in 3D and appeal to the celebrated two-tier energy method developed by Guo and Tice \cite{GuTi1,GuTi2} in the context of viscous surface waves without surface tension. The conclusion is that we are able to prove the global well-posedness for 3D full compressible viscous non-resistive MHD system with strong background magnetic field and furthermore show that the unique solution converges at an almost exponential rate to the steady state. In particular, our result reveals in some sense that the vertical magnetic field enjoys the stabilization effects for the compressible viscous non-resistive fluids, even if the temperature equation is incorporated.
It should be stressed that the two-tier energy method is quite effective in studying the Rayleigh-Taylor instability of compressible/incompressible MHD system. We refer to a series of work by Jiang and Jiang \cite{JJ1,JJ2,JJ3,JJ4} and Wang \cite{Wa1} for more results in this direction.

\subsection{Reformulations in Lagrangian coordinates}\label{re_l}
In analogy with \cite{LXZ,XZ,TW}, it is more convenient and effective to work with Lagrangian coordinates. To do this, we first assume that there exists an invertible mapping $\eta_0:\Om \rightarrow \Om$ such that
$\p \Om =\eta_0(\p \Om)$. The flow mapping $\eta$ is then defined through
\begin{equation}\label{re_l-1}
\left\{\begin{aligned}
& \p_t \eta(t,x)=\tilde{\vu}(t,\eta(t,x)), \\
& \eta|_{t=0}=\eta_0.\\
\end{aligned}\right.
\end{equation}

Obviously, if the flow mapping $\eta$ is a small perturbation of the identity mapping $Id$, the mapping $\eta(t,\cdot)$ is invertible and a diffeomorphism for any $t\in \R^{+}$.\footnote{It then follows from (\ref{fmhd_7})$_1$, (\ref{re_l-1}) and the hypothesis $\p \Om =\eta_0(\p \Om)$ that $\p \Om=\eta(t,\p \Om)$.} Consequently, the unknowns in Lagrangian coordinates read as
\begin{equation}\label{re_l-2}
(\vr,\vu,\vb,\vt)(t,x)=(\tilde{\vr},\tilde{\vu},\tilde{\vb},\tilde{\vt})(t,\eta(t,x)),\,\, t\in \R^{+},x\in \Om.
\end{equation}

In terms of $(\eta,\vr,\vu,\vb,\vt)$, the equations (\ref{fmhd_1})-(\ref{fmhd_5}) are thus reformulated as
\begin{equation}\label{re_l-3}
\p_t \eta=\vu,
\end{equation}
\begin{equation}\label{re_l-4}
\p_t \vr+\vr \Div_{\mathcal{A}} \vu=0,
\end{equation}
\begin{equation}\label{re_l-5}
\vr \p_t \vu -\mu \Delta_{\mathcal{A}}\vu-(\mu+\lambda)\Grad_{\mathcal{A}}\Div_{\mathcal{A}}\vu +\Grad_{\ma}\left( \vr\vt+\f{1}{2}|\vb|^2 \right)=\vb\cdot \Grad_{\ma}\vb,
\end{equation}
\begin{equation}\label{re_l-6}
\p_t \vb-\vb\cdot \Grad_{\ma}\vu +\vb \Div_{\mathcal{A}} \vu=\mathbf{0},
\end{equation}
\begin{equation}\label{re_l-7}
\vr \p_t \vt+\vr\vt \Div_{\mathcal{A}} \vu
=\kappa \Delta_{\ma} \vt+
2\mu |\mathbb{D}_{\ma}\vu|^2+\lambda(\Div_{\mathcal{A}} \vu)^2,
\end{equation}
\begin{equation}\label{re_l-8}
\Div_{\mathcal{A}} \vb=0,
\end{equation}
in $\R^{+}\times \Om$. It is clear that the initial-boundary conditions reduce to
\begin{equation}\label{re_l-9}
(\eta,\vr,\vu,\vb,\vt)|_{t=0}=(\eta_0,\vr_0,\vu_0,\vb_0,\vt_0)
\eeq
with $(\vr_0,\vu_0,\vb_0,\vt_0)(x):=(\tilde{\vr}_0,\tilde{\vu}_0,\tilde{\vb}_0,\tilde{\vt}_0)(\eta_0(x)),x\in \Om$, as well as
\begin{equation}\label{re_l-10}
\left\{\begin{aligned}
& \vu|_{\p \Om}=\mathbf{0}, \\
& \vt|_{\p \Om}=\overline{\vt}. \\
\end{aligned}\right.
\end{equation}

In the formulation above we denote by $\ma=( (\Grad_{x}\eta)^{-1}  )^{t}=:(\ma_{ij})_{i,j=1}^3$ and the differential operators subscripted with $\ma$ means the actions with respect to $\ma$. More precisely,
\[
(\Grad_{\ma}g)_{i}:=\ma_{ij}\p_j g,\,\, \,\, \Div_{\ma}\mathbf{G}:=\ma_{ij}\p_j G_i , \text{  with   } \mathbf{G}=(G_1,G_2,G_3)^{t},
\]
\[
\Delta_{\ma} g = \Div_{\ma}\Grad_{\ma}g,\,\, \,\,(\mathbb{D}_{\ma}\vu)_{ij}=\f{   \ma_{ik}\p_k u_j+ \ma_{jk}\p_k u_i       }{2} .
\]

Let $\overline{\vr}>0, \overline{\vt}>0,\overline{b} \in \R, \overline{b}\neq 0$ be constants and $\overline{\vb}:=(0,0,\overline{b})^{t}$. We shall prove global well-posedness of the initial-boundary value problem (\ref{re_l-3})-(\ref{re_l-10}) as long as the initial data satisfy suitable compatibility conditions specified later in the main theorem, sufficiently smooth and close enough to the steady state $(Id,\overline{\vr},\mathbf{0},\overline{\vb},\overline{\vt})$. Moreover, the unique global solution converges to the steady state at an almost exponential rate as $t\rightarrow \infty$.

Following \cite{LXZ,TW,XZ}, we furthermore make use of several conserved quantities in order to capture the weak dissipation of the magnetic field. The interested reader may consult \cite{TW} for the details. We denote by $J:=\text{det}( \Grad_x \eta)$ the Jacobian of the coordinate transformation. A straightforward computation gives
\beq\label{re_l-11}
\p_t J= J \Div_{\ma}\vu.
\eeq
Combining (\ref{re_l-11}) and the equation of continuity (\ref{re_l-4}) we arrive at
\beq\label{re_l-12}
\p_t (\vr J)=0.
\eeq
In addition, testing the equation of magnetic field (\ref{re_l-6}) by $J\ma^{t}$ yields
\beq\label{re_l-13}
\p_t ( J\ma^{t} \vb )=\mathbf{0}.
\eeq
It is clear that
\beq\label{re_l-14}
\p_t \eta|_{\p \Om}=\mathbf{0}.
\eeq
Bearing in mind that we expect the solution converges to the steady state $(Id,\overline{\vr},\mathbf{0},\overline{\vb},\overline{\vt})$ as time tends to infinity, we may deduce from (\ref{re_l-12})-(\ref{re_l-14}) that
\beq\label{re_l-15}
\vr J= \overline{\vr},\,\,J \ma^{t}\vb = \overline{\vb} \,\, \text{  in  }\Om,\,\, \,\, \,\, \,\, \,\,  \eta =Id  \,\, \text{  on  }\p \Om.
\eeq
Notice also that there holds $\Div_{\ma}\vb=J^{-1}\Div_x(J\ma^{t}\vb)=J^{-1}\Div_x \overline{\vb}=0$; whence the divergence-free condition of magnetic field (\ref{re_l-8}) holds. For convenience, we make a shift of the flow mapping via $\eta \rightarrow \eta+Id $ and recover the relations $J=\text{det}( I+ \Grad_x \eta)$, $\ma=( (I+\Grad_{x}\eta)^{-1}  )^{t}$. Obviously, (\ref{re_l-15}) is reformulated as
\beq\label{re_l-16}
\vr=\overline{\vr}J^{-1},\,\, \vb=J^{-1}(I+ \Grad_x \eta)\overline{\vb}=\overline{b}J^{-1}(\mathbf{e}_3 +\p_3 \eta   ) \,\, \text{  in  }\Om,\,\, \,\, \,\, \,\,
\eta =\mathbf{0}  \,\, \text{  on  }\p \Om.
\eeq
Based on (\ref{re_l-16}), one finds after a direct calculation that
\beq\label{re_l-17}
\vb\cdot \Grad_{\ma} \vb- \Grad_{\ma}\left(  \f{|\vb|^2}{2}   \right)=
\overline{b}^2 (   \p_{33}\eta -\p_3 \Div_x\eta \mathbf{e}_3 +\Grad_x \Div_x \eta- \Grad_x \p_3 \eta_3     ) +\mathbf{R}^{\eta},
\eeq
with $\mathbf{R}^{\eta}=O( \Grad_x \eta \Grad_x^2 \eta )$ being the remainder. Finally, we also take a shift of the temperature via $\vt\rightarrow \vt+\overline{\vt}$. As a consequence, we reformulate (\ref{re_l-3})-(\ref{re_l-10}) in $\R^{+}\times \Om$ as
\begin{equation}\label{re_l-18}
\left\{\begin{aligned}
& \p_t \eta=\vu, \\
& \overline{\vr}J^{-1}\p_t \vu -\mu \Delta_{\ma}\vu -(\mu+\lambda)\Grad_{\ma}\Div_{\ma}\vu +   \Grad_{\ma} [  \overline{\vr}J^{-1}  (\vt+\overline{\vt}) ]        \\
&= \overline{b}^2 (   \p_{33}\eta -\p_3 \Div_x\eta \mathbf{e}_3 +\Grad_x \Div_x \eta- \Grad_x \p_3 \eta_3     ) +\mathbf{R}^{\eta},\\
& \overline{\vr}J^{-1} \p_t \vt+   \overline{\vr}J^{-1}(\vt+\overline{\vt}) \Div_{\mathcal{A}} \vu
=\kappa \Delta_{\ma} \vt+
2\mu |\mathbb{D}_{\ma}\vu|^2+\lambda(\Div_{\mathcal{A}} \vu)^2,\\
& (\eta,\vu,\vt)|_{t=0}=(\eta_0,\vu_0,\vt_0),\\
& \vu|_{\p \Om}=\mathbf{0},\,\,\,  \vt |_{\p \Om}=0. \\
\end{aligned}\right.
\end{equation}

The rest of this paper is devoted to showing the global well-posedness for (\ref{re_l-18}) around the steady state $(\mathbf{0},\mathbf{0},0)$ and that the solution converges at an almost exponential rate. Moreover, with the solution $(\eta,\vu,\vt)$ of (\ref{re_l-18}) at hand, we may define $\vr=\overline{\vr}\text{det}(I+\Grad_x \eta)^{-1},\,\vb=\overline{b}\text{det}(I+\Grad_x \eta)^{-1}(\mathbf{e}_3 +\p_3 \eta   )$ and conclude that $(\eta+Id,\vr,\vu,\vb,\vt+\overline{\vt})$ solves (\ref{re_l-3})-(\ref{re_l-10}) with initial data $(\eta_0+Id,\vr_0,\vu_0,\vb_0,\vt_0+\overline{\vt})$, where $\vr_0=\overline{\vr}\text{det}(I+\Grad_x \eta_0)^{-1},\,\vb_0=\overline{b}\text{det}(I+\Grad_x \eta_0)^{-1}(\mathbf{e}_3 +\p_3 \eta_0   )$.

\section{Main theorem and strategy of the proof}\label{ma_sta}

We first introduce the notations used from now on. Let $1\leq p\leq \infty,k,l\in \mathbb{N}$. We denote by $L^p(\Om)$ and $H^k(\Om)$ the usual Lebesgue space and Sobolev space respectively; the norms are denoted by $\|\cdot \|_{L^p}$ and $\| \cdot \|_{k}$ respectively. Moreover, for convenience we define the following anisotropic Sobolev norm
\[
\| g \|_{k,l}:= \sum_{\alpha_1+\alpha_2\leq l} \| \p _1^{\alpha_1}  \p_2^{\alpha_2}g \|_{k} \,\,\,\,   \text{   for appropriate } g.
\]
In particular, $\| \cdot \|_{k,0}=\|\cdot \|_{k}$. $\Grad_{\ast},\Div_{\ast},\Delta_{\ast}$ means the operators acting only on the horizontal variables, i.e.,
\[
\Grad_{\ast}g=(\p_1 g, \p_2 g)^{t},\,\,\,\, \Delta_{\ast} g =\p_{11} g+ \p_{22} g,\,\,\,\,
\Div_{\ast} \mathbf{G}=\p_1 G_1 +\p_2 G_2 \text{  with   } \mathbf{G}=(G_1,G_2,G_3)^{t}.
\]
For vector $\mathbf{G}=(G_1,G_2,G_3)^{t}$, we denote by $\mathbf{G}_{\ast}:=(G_1,G_2)^{t}$ the horizontal components. Let $l\in \{ 2,3\}$. The set of time-space multi-indices is denoted by
$\mathbb{N}^{1+l}:= \{     \alpha=(\alpha_0,\alpha_1,...,\alpha_{l})   \}        $ with parabolic counting $|\alpha|=2\alpha_0+\alpha_1+...+\alpha_{l}$. Similarly, the set of space multi-indices is denoted by $\mathbb{N}^{l}:=\{
 (\alpha_1,...,\alpha_{l})  \}$ with $|\alpha|=\alpha_1+...+\alpha_{l}$. For convenience, we set
 \[
 \| \overline{\Grad}^k g \|_0^2 := \sum_{     \alpha \in \mathbb{N}^{1+3},|\alpha|\leq k        }  \|      \p_t^{\alpha_0}\p_1^{\alpha_1}\p_2^{\alpha_2}\p_3^{\alpha_3} g\|_0^2.
 \]
We denote by $C$ a generic positive constant independent of the initial data and time, changing from line to line. For brevity, we use $\int_{\Om} f$ for the integration of $f$ over $\Om$, while neglecting the usual $\dx$.
We frequently use $A\lesssim B$ to denote $A\leq C B$ for generic positive constant $C$.

\subsection{Main theorem}\label{ma_th}

Before stating the main result, we introduce the energy functionals. Let $n =2N$ or $N+2$ with $N \geq 4$ being an integer. We set the energy as
\beq\label{ma_th-1}
\mathcal{E}_n:= \sum_{j=0}^n\| \p_t ^j \vu \|_{2n-2j}^2   +\sum_{j=0}^n\| \p_t ^j \vt \|_{2n-2j}^2
+\| \eta \|^2_{2n+1},
\eeq
and the dissipation as
\beq\label{ma_th-2}
\mathcal{D}_n:= \sum_{j=0}^n\| \p_t ^j \vu \|_{2n-2j+1}^2  +\sum_{j=0}^n\| \p_t ^j \vt \|_{2n-2j+1}^2
+\| \Div_x \eta \|^2_{2n}+ \| \p_3 \eta \|^2_{2n}+\|  \eta \|^2_{2n}.
\eeq
Finally, we set
\begin{equation}\label{ma_th-3}
     \begin{split}
       \mathcal{G}_{2N}(t):= \sup_{0\leq \tau \leq t}\mathcal{E}_{2N}(\tau)  & \       +\int_0^{t}\mathcal{D}_{2N}(\tau) \rm{d} \tau         \\
    & \   +\sup_{0\leq \tau \leq t}  (1+\tau)^{2N-4} \mathcal{E}_{N+2}(\tau)  \\
     \end{split}
\end{equation}

Our global well-posedness and decay of solutions to (\ref{re_l-18}) read as
\begin{Theorem}\label{m-th}
Let $N \geq 4$ be an integer. Suppose that $\vu_0\in H^{4N}(\Om),\,\vt_0\in  H^{4N}(\Om),\, \eta_0 \in H^{4N+1}(\Om),\,\eta_0|_{\p \Om}=\mathbf{0}$ satisfy the following natural compatibility conditions:
\begin{equation}\label{ma_th-3-1}
\p_t^j \vu (0,x)=\mathbf{0},\,\,\,\,\p_t^j \vt(0,x)=0,\,\,\,\,\text{         for any      }  x\in \p \Om,
\end{equation}
with $j=0,1,...,2N-1$. Then there exists $\ep_0>0$ such that if $\mathcal{E}_{2N}(0)\leq \ep_0$ the problem (\ref{re_l-18}) admits a unique global solution $(\eta,\vu,\vt)$ on $[0,\infty)$. Moreover, it holds the uniform estimate
\beq\label{ma_th-4}
\mathcal{G}_{2N}(\infty) \lesssim \mathcal{E}_{2N}(0).
\eeq

\end{Theorem}

Several remarks concerning the main result are as follows.

\begin{Remark}
To simplify the presentation we have chosen the background magnetic field $\overline{\vb}$ to be $(0,0,\overline{b})^{t}$ with $\overline{b}\neq 0$. As a matter of fact, one checks by following the same line of the proof below
that the general case $\overline{\vb}=(\overline{B}_1,\overline{B}_2,\overline{B}_3)^{t}$ with $\overline{B}_3 \neq 0$ could be treated exactly in the same manner.

\end{Remark}

\begin{Remark}
We notice from the uniform estimate (\ref{ma_th-4}) that $\mathcal{E}_{N+2}(t) \lesssim (1+t)^{-2N+4}$, which implies the almost exponential decay estimate of solution since $N$ can be any integer greater than $4$, as long as the initial data are sufficiently smooth.
In particular, $\| \eta(t)\|_{2N+5}$ is sufficiently small for any $t\geq 0$, ensuring the flow mapping a diffeomorphism. Thus, Theorem \ref{m-th} yields the global well-posedness of (\ref{fmhd_1})-(\ref{fmhd_7}) in Eulerian coordinates.

\end{Remark}

\begin{Remark} \label{re_3}
We remark that the local well-posedness to the full compressible viscous non-resistive MHD system is nowadays standard. On the one hand, we could directly prove the local well-posedness to the geometric form (\ref{re_l-18}) based on linearization, construction of approximating solutions and application of fixed point argument. On the other hand, recalling that we already have the local well-posedness in Eulerian coordinates \cite{Kaw1,FnYu1}, the corresponding result in Lagrangian coordinates follows immediately since the smallness of solutions ensures that the flow mapping $\eta(t,\cdot):\Om \rightarrow \Om$ is a diffeomorphism for any $t\in \R^{+}$. Therefore, we shall omit the details of proof for local well-posedness, whereas focus attention on the global energy estimates in the sequel.
\end{Remark}

\subsection{Strategy of the proof}\label{st_pr}

Basically, the strategy of the proof of Theorem \ref{m-th} relies on the celebrated two-tier energy method of Guo, Tice \cite{GuTi1,GuTi2} and follows the same line of Tan, Wang \cite{TW} for barotropic fluids with necessary modifications so as to accommodate the appearance of temperature. More precisely,
\begin{enumerate}[(1)]
\item {  We linearize the system (\ref{re_l-18}) in the geometric form and derive the energy evolution of time derivatives, see Section \ref{e-e-t}.                     }
\item {  We reformulate the system (\ref{re_l-18}) in the linear perturbed form and obtain the energy evolution of horizontal spatial derivatives, see Section \ref{e-e-h}. Notice also that the temporal derivatives and horizontal derivatives naturally preserve the homogeneous boundary conditions; whence it is convenient to integrate by parts and apply the Poincar\'{e} inequality.         }
\item {  In order to obtain the dissipation estimates of the flow mapping, we take $\p ^{\alpha} \eta $ to be the test function in the momentum equation and make use of the structure of equations, see Section \ref{e-e-f}.             }
\item{ By exploiting the structure between the first two components and the third component of the momentum equation, we are able to derive the energy evolution of vertical derivative of solutions, see Section \ref{e-e-v}.
    }
\item{ We improve the estimates for $(\p_t \vu,\p_t \vt)$ in high Sobolev norms via the classical elliptic regularity theory in a suitable manner, see Section \ref{e-i-e}.
}
\item{ We glue all the energy evolution estimates at hand to show first the boundedness of $\mathcal{E}_{2N}$ and then the decay estimates of $\mathcal{E}_{N+2}$. This closes the global energy estimates under the assumption that $\mathcal{G}_{2N}(T)\leq \delta$ for $\delta>0$ sufficiently small, see Section \ref{g-e}. This completes the proof of our main theorem.
}
\end{enumerate}

Several crucial observations of our paper read as follows.
\begin{itemize}
\item{
To linearize the pressure in the geometric form (see (\ref{e-e-t-2})), we adopt the elementary identities
\begin{equation}\nonumber
     \begin{split}
         \Grad_{\ma}[ \overline{\vr} J^{-1} (\vt+\overline{\vt})  ] & \         =   \Grad_{\ma}   [ \overline{\vr} J^{-1} \vt- \overline{\vr}\vt+
              \overline{\vt} (  \overline{\vr} J^{-1} - \overline{\vr}) +  \overline{\vr} (\vt+\overline{\vt}  )          ]               \\
    & \    =   \Grad_{\ma}    [  (  \overline{\vr} J^{-1} - \overline{\vr} )\vt                          ]  + \overline{\vt}  \Grad_{\ma}  (  \overline{\vr} J^{-1} - \overline{\vr})
    +    \overline{\vr}   \Grad_{\ma}  \vt    \\
         & \  =   \overline{\vr}   \Grad_{\ma}  \vt -  \overline{\vr}\overline{\vt} \Grad_x \Div_x \eta+ \Grad_{\ma}    [  (  \overline{\vr} J^{-1} - \overline{\vr} )\vt                          ] + \mathbf{R}^{\eta} ,              \\
     \end{split}
\end{equation}
with $\mathbf{R}^{\eta}=O( \Grad_x \eta \Grad_x^2 \eta )$ being the remainder. The first two linear terms on the right-hand side are thus taken to be the linearized pressure. It should be emphasized that the term $-  \overline{\vr}\overline{\vt} \Grad_x \Div_x \eta$ appears naturally in our framework of Lagrangian coordinates and plays a key role in the analysis, see for instance Lemma \ref{le-7}.
}
\item{
We make use of the cancellations between the linearized momentum equation and the temperature equation when deriving the energy evolution of time derivatives and the horizontal derivatives, cf. (\ref{e-e-t-14}).
}
\item{
When deriving the energy evolution of the flow mapping $\eta$, see Lemma \ref{le-6}, we encounter an extra integral $\int_{\Om} \overline{\vr}  \Grad_x (\p^{\alpha}\vt)\cdot \p^{\alpha}\eta $ due to the effect of temperature. However, such an integral cannot be controlled in the same way as the nonlinear perturbation term $\int_{\Om} \p^{\alpha}\eta \cdot \p^{\alpha}\mathbf{W}$. To this end, we split it into two cases, according to either $|\alpha|=0$ or $0<|\alpha|\leq 2n$. In the first case, we estimate through the recovering dissipation $ \mathcal{D}_{n}^{ (\sharp) }$ and the temporal dissipation $\mathcal{D}_{n}^{ (t) }$ suitably; in the second case, we estimate through the recovering dissipation $ \mathcal{D}_{n}^{ (\sharp) }$ and the horizontal dissipation $\mathcal{D}_{n}^{ (\ast) }$ suitably, cf. (\ref{e-e-f-10})-(\ref{e-e-f-11}).
}

\item{
Due to the strong coupling between the velocity and temperature fields, the additional term $ \sum_{j=1}^{n-1} \|  \p_t^j \vt \|  ^2_{2n-2j}    $ appears when one improves the estimate $ \sum_{j=1}^{n} \|  \p_t^j \vu \|  ^2_{2n-2j+1}    $; this new term cannot be controlled as in the barotropic case, cf. Proposition 3.8 in \cite{TW}. To this end, we make use of the inner structure of the system. More precisely, we improve the estimate $ \sum_{j=0}^{n} \|  \p_t^j \vt \|  ^2_{2n-2j+1}    $ simultaneously. When controlling the latter we also encounter the term $ \sum_{j=1}^{n-1}    \| \p_t ^j \vu  \|^2_{2n-2j} $ involving the velocity field. We finally close this process with the help of Sobolev interpolation inequality and Young inequality, see Lemma \ref{le-8} for more details. Remarkably, the same spirit is applied also in Lemma \ref{le-9}.
}

\end{itemize}

\section{Energy evolution}\label{e-e}
The key issue of proving Theorem \ref{m-th} is to derive sufficient global a priori estimates. To this end, we assume that $(\eta,\vu,\vt)$ is the smooth solution of problem (\ref{re_l-18}) and the solution satisfies the estimate $\mathcal{G}_{2N}(T)\leq \delta$ for some $\delta>0$ suitably small.

\subsection{Energy evolution of time derivatives}\label{e-e-t}
In this subsection, we derive the estimates for time derivatives of solutions in the geometric form. To do that we first notice the basic identity
\begin{equation}\label{e-e-t-1}
J^{-1}=1-\Div_x \eta +O(|\Grad_x \eta|^2)
\end{equation}
and linearize (\ref{re_l-18}) as follows
\begin{equation}\label{e-e-t-2}
\left\{\begin{aligned}
& \p_t \eta=\vu, \\
& \overline{\vr}J^{-1}\p_t \vu -\mu \Delta_{\ma}\vu -(\mu+\lambda)\Grad_{\ma}\Div_{\ma}\vu +  \overline{\vr} \Grad_{\ma} \vt  -  \overline{\vr}  \overline{\vt} \Grad_x \Div_x \eta   \\
&= \overline{b}^2 (   \p_{33}\eta -\p_3 \Div_x\eta \mathbf{e}_3 +\Grad_x \Div_x \eta- \Grad_x \p_3 \eta_3     ) +\mathbf{R}^{\eta}+\mathbf{F},\\
& \overline{\vr}J^{-1} \p_t \vt+   \overline{\vr}\overline{\vt} \Div_{\mathcal{A}} \vu
=\kappa \Delta_{\ma} \vt+G ,\\
& (\eta,\vu,\vt)|_{t=0}=(\eta_0,\vu_0,\vt_0),\\
& \vu|_{\p \Om}=\mathbf{0},\,\,\,  \vt |_{\p \Om}=0. \\
\end{aligned}\right.
\end{equation}
with
\begin{equation}\label{e-e-t-3}
\mathbf{F}:=-\Grad_{\ma}  [    \overline{\vr} (J^{-1}-1) \vt  ],
\end{equation}
\begin{equation}\label{e-e-t-4}
G:=-  \Big(
\overline{\vr} \vt+ \overline{\vr} (J^{-1}-1)\vt + \overline{\vr} (J^{-1}-1) \overline{\vt}
\Big) \Div_{\mathcal{A}} \vu +2\mu |\mathbb{D}_{\ma}\vu|^2+\lambda(\Div_{\mathcal{A}} \vu)^2.
\end{equation}

For $j=0,...,n$, we apply the time derivative $\p_t^{j}$ to the linearized system (\ref{e-e-t-2}) to conclude that
\begin{equation}\label{e-e-t-5}
\left\{\begin{aligned}
& \p_t (\p_t^{j}\eta)=\p_t^{j}\vu, \\
& \overline{\vr}J^{-1}\p_t (\p_t^{j}\vu) -\mu \Delta_{\ma}(\p_t^{j}\vu) -(\mu+\lambda)\Grad_{\ma}\Div_{\ma}(\p_t^{j}\vu) +  \overline{\vr} \Grad_{\ma} (\p_t^{j}\vt)  -  \overline{\vr}  \overline{\vt} \Grad_x \Div_x (\p_t^{j}\eta )  \\
&= \overline{b}^2 \Big(   \p_{33}(\p_t^{j}\eta) -\p_3 \Div_x  (\p_t^{j}\eta) \mathbf{e}_3 +\Grad_x \Div_x (\p_t^{j}\eta)- \Grad_x \p_3 (\p_t^{j}\eta_3)     \Big) +\mathbf{N}^j,\\
& \overline{\vr}J^{-1} \p_t (\p_t^{j}\vt)+   \overline{\vr}\overline{\vt} \Div_{\mathcal{A}} (\p_t^{j}\vu)
=\kappa \Delta_{\ma}  (\p_t^{j}\vt)+M^j ,\\
& \p_t^{j}\vu|_{\p \Om}=\mathbf{0},\,\,\,  \p_t^{j}\vt |_{\p \Om}=0. \\
\end{aligned}\right.
\end{equation}
Here the nonlinear term $\mathbf{N}^j=(N^j_1,N^j_2,N^j_3)$ with
\begin{equation}\label{e-e-t-6}
     \begin{split}
         N^j_i=\sum_{0< \ell \leq j} & \   C_j^{\ell} \Big\{  \mu \ma_{lk} \p_k ( \p_t^{\ell}\ma_{lm}\p_t ^{j-\ell}\p_m u_i    ) +\mu \p_t^{\ell} \ma_{lk} \p_t^{j-\ell}\p_k (\ma_{lm}\p_m u_i) \\
     & \   +(\mu+\lambda)\ma_{ik}\p_k ( \p_t^{\ell} \ma_{lm}\p_t^{j-\ell}\p_m u_l   )+(\mu+\lambda) \p_t^{\ell} \ma_{ik} \p_t^{j-\ell} \p_k ( \ma_{lm}\p_m u_l  )          \\
         & \   -\overline{\vr}\p_t^{\ell}(J^{-1}) \p_t (\p_t^{j-\ell}u_i )  -\overline{\vr} \p_t^{\ell} \ma_{ik} \p_t^{j-\ell}\p_k \vt \Big\}           \\
         & \     - \sum_{0\leq \ell \leq j} C_j^{\ell}\p_t ^{\ell} \ma_{ik} \p_t ^{j-\ell} \p_k[    \overline{\vr} (J^{-1}-1) \vt  ]        \\
         & \  +\p_t^j R^{\eta}_i, \,\,\,\,\,\, i=1,2,3,
     \end{split}
\end{equation}
as well as
\begin{equation}\label{e-e-t-7}
     \begin{split}
         M^j=\sum_{0< \ell \leq j} & \    C_j^{\ell}  \Big\{  \mu \ma_{lk} \p_k ( \p_t^{\ell}\ma_{lm}\p_t ^{j-\ell}\p_m \vt    ) +\mu \p_t^{\ell} \ma_{lk} \p_t^{j-\ell}\p_k (\ma_{lm}\p_m \vt) \\
          & \    -   \overline{\vr}\p_t^{\ell}(J^{-1}) \p_t (\p_t^{j-\ell}\vt)- \overline{\vr} \overline{\vt}\p_t^{\ell}\ma_{ik}\p_t ^{j-\ell}\p_k u_i \Big\}            \\
            &\  -   \sum_{0\leq \ell \leq j} C_j^{\ell}  \Big\{   \overline{\vr} \p_t^{\ell} \vt \p_t^{j-\ell} (\ma_{ik} \p_k u_i)+
            \overline{\vr}     \p_t^{\ell} [ (J^{-1}-1)\vt] \p_t^{j-\ell} (\ma_{ik} \p_k u_i)                \\
            &\   +  \overline{\vr}  \overline{\vt}  \p_t^{\ell} (J^{-1}-1) \p_t^{j-\ell} (\ma_{ik} \p_k u_i)\Big\} \\
            &\  +\p_t^j   \Big\{ \f{\mu}{2}  |     \ma_{ik}\p_k u_l +\ma_{lk}\p_k u_i               |^2 +\lambda (\ma_{ik}\p_k u_i)^2 \Big\}           .
     \end{split}
\end{equation}

We have the following estimates on the nonlinear terms $\mathbf{N}^j$ and $M^j$; these are quite useful for future analysis.
\begin{Lemma}\label{le-1}
Let $n\geq 3$ be an integer. It holds
\begin{equation}\label{e-e-t-8}
\|  \mathbf{N}^j  \|_{0}^2 +\| M^j \|_0^2
\lesssim \mathcal{E}_n \mathcal{D}_n.
\end{equation}
\end{Lemma}
{\bf{Proof.} } Expanding $N^j_i$ further by Leibniz rule, we find that each term is at least quadratic and could be written as $\mathscr{X}\mathscr{Y}$, where $\mathscr{X}$ contains fewer derivatives than $\mathscr{Y}$.
We use the Sobolev embeddings and the definitions of $\mathcal{E}_n,\mathcal{D}_n$ to obtain
\[
\| \mathscr{X} \|_{L^{\infty}} \lesssim \sqrt{  \mathcal{E}_n},\,\,\,\, \| \mathscr{Y} \|_{0} \lesssim \sqrt{ \mathcal{D}_n};
\]
whence
\[
\| \mathscr{X}\mathscr{Y} \|_0^2 \leq \| \mathscr{X} \|_{L^{\infty}}^2 \| \mathscr{Y} \|_{0} ^2  \lesssim \mathcal{E}_n \mathcal{D}_n.
\]
The estimate of $\| M^j \|_0^2$ is carried out in the same way. This gives (\ref{e-e-t-8}).    $\Box$

We introduce the temporal energy as
\begin{equation}\label{e-e-t-9}
     \begin{split}
       \mathcal{E}_n^{(t)}=\sum_{j=0}^n & \   \Big\{   \overline{\vr} \|\p_t^j \vu \|_0^2 +\overline{\vr} \overline{\vt}^{-1} \|\p_t^j \vt \|_0^2
       +  \overline{\vr}  \overline{\vt} \| \p_t^j \Div_x \eta   \|_0^2    \\
     & \   +\overline{b}^2 \| \p_t^j \p_3 \eta_{\ast} \|_0^2 +\overline{b}^2 \| \p_t^j \Div_{\ast} \eta_{\ast} \|_0^2 \Big\},  \\
     \end{split}
\end{equation}
together with the temporal dissipation as
\begin{equation}\label{e-e-t-10}
\mathcal{D}_n^{(t)} =\sum_{j=0}^n  \left(
\| \p_t^j \vu \|_1^2+\| \p_t^j \vt \|_1^2
\right).
\end{equation}
Now we are ready to give the energy evolution of temporal derivatives. Compared with the barotropic case \cite{TW}, we need some cancellations between the momentum equation and the temperature equation due to the strong coupling effects between the velocity and temperature.
\begin{Lemma}\label{le-2}
Let $n\geq 3$ be an integer. Then
\begin{equation}\label{e-e-t-11}
\f{\rm{d}}{\dt} \mathcal{E}_n^{(t)} +\mathcal{D}_n^{(t)} \lesssim
\sqrt{\mathcal{E}_n}  \mathcal{D}_n.
\end{equation}
\end{Lemma}
{\bf{Proof.} } Taking $L^2$ inner products of (\ref{e-e-t-5})$_2$ with $J \p_t^j \vu$ and (\ref{e-e-t-5})$_3$ with $\overline{\vt}^{-1}J\p_t^j \vt$ respectively, followed by adding the resulting relations, we obtain after integration by parts that
\begin{equation}  \label{e-e-t-12}
     \begin{split}
         \f{1}{2} \f{\rm{d}}{\dt} \int_{\Om} \left( \overline{\vr}|\p_t^j \vu|^2+\overline{\vr} \overline{\vt}^{-1}  |\p_t^j \vt|^2\right) & \  + \int_{\Om}
         \Big\{  \mu J |\Grad_{\ma}\p_t^j \vu|^2  +(\mu+\lambda)J|\Div_{\ma}\p_t^j \vu|^2    \\
         &\ +\kappa \overline{\vt}^{-1} J |\Grad_{\ma}\p_t^j \vt|^2  \Big\} \\
         & \   +\int_{\Om} \overline{\vr}\overline{\vt}\Div_x(\p_t^j \eta) \Div_x(\p_t^j \vu)      \\
         & \   +\int_{\Om}\overline{b}^2  \Big\{  \p_3 (\p_t^j \eta)\cdot \p_3(\p_t^j \vu)-\Div_x(\p_t^j \eta) \p_3 (\p_t^j u_3)        \\
         &  \   +  \Div_x(\p_t^j \eta) \Div_x(\p_t^j \vu)   -\p_3 (\p_t^j \eta_3) \Div_x(\p_t^j \vu) \Big\}           \\
         = & \ \int_{\Om} \left( J \p_t^j \vu \cdot \mathbf{N}^j+ \overline{\vt}^{-1} J \p_t^j \vt M^j\right)         \\
         & \ -\int_{\Om}  \overline{\vr}\overline{\vt}   \Div_x(\p_t^j \eta) \Div_x ((J-1)\p_t^j \vu  )  \\
         & \ -\int_{\Om}  \overline{b}^2 \Big\{  \p_3(\p_t^j \eta) \cdot \p_3 ( (J-1)\p_t^j \vu  ) - \Div_x(\p_t^j \eta) \p_3 ( (J-1)\p_t^j u_3      ) \\
         & \  + \Div_x(\p_t^j \eta) \Div_x (  (J-1)\p_t^j \vu     ) -\p_3(  \p_t^j \eta_3)\Div_x (  (J-1)\p_t^j \vu     )  \Big\}.
     \end{split}
\end{equation}
Indeed,
\begin{equation}\label{e-e-t-13}
     \begin{split}
         -\int_{\Om} \Delta_{\ma} (\p_t^j \vt) J \p_t^j \vt &\  =-\int_{\Om} \Div_{\ma}(\Grad_{\ma} (\p_t^j \vt)) J \p_t^j \vt  \\
     & \  =-   \int_{\Om} \ma_{il}\p_l (\ma_{ik}\p_k(\p_t^j \vt) )  J\p_t^j \vt    \\
         & \ =   \int_{\Om} \ma_{ik}\p_k(\p_t^j \vt) \p_l(J\p_t^j \vt \ma_{il} )        \\
         & \  =\int_{\Om} \ma_{ik}\p_k(\p_t^j \vt) \Big( J \ma_{il}  \p_l (\p_t^j \vt)+\p_t^j \vt \p_l(J\ma_{il})   \Big) \\
         & \ =\int_{\Om} J |\Grad_{\ma}\p_t^j \vt|^2,
     \end{split}
\end{equation}
due to the facts that $\p_t^j \vt|_{\p \Om}=0$ and $\p_k(J\ma_{ik})=0$; similarly,

\begin{equation}\label{e-e-t-14}
     \begin{split}
        \int_{\Om}\overline{\vr} \Grad_{\ma}  ( \p_t^j \vt)  J \p_t^j \vu & \  = \int_{\Om}\overline{\vr}
        \ma_{ik}\p_k(\p_t^j \vt) J \p_t^j u_i       \\
              & \  =- \int_{\Om}\overline{\vr}  \p_t^j \vt \p_k ( J \p_t^j u_i   \ma_{ik}     )             \\
         & \   =- \int_{\Om}\overline{\vr}  \p_t^j \vt  \Big( \p_t^j u_i   \p_k(  J  \ma_{ik}  ) +  J  \ma_{ik}  \p_k (\p_t^j u_i)      \Big)  \\
         & \  =-\int_{\Om}\overline{\vr}  J \Div_{\ma}(\p_t^j \vu) \p_t^j \vt,        \\
     \end{split}
\end{equation}
which cancels with the term $\overline{\vr}\overline{\vt} \Div_{\ma}(\p_t^j \vu) $ in (\ref{e-e-t-5})$_3$ upon multiplying by $\overline{\vt}^{-1}J\p_t^j \vt$ and integration over $\Om$.

Next, we aim to simplify the expression of (\ref{e-e-t-12}). Using (\ref{e-e-t-5})$_1$,
\begin{equation}\label{e-e-t-15}
     \begin{split}
         \int_{\Om} \overline{\vr}\overline{\vt}\Div_x(\p_t^j \eta) \Div_x(\p_t^j \vu)  & \ = \int_{\Om} \overline{\vr}\overline{\vt}\Div_x(\p_t^j \eta) \Div_x(\p_t^{j+1} \eta)  \\
     & \ =\overline{\vr}\overline{\vt}\f{1}{2} \f{\rm{d}}{\dt}\int_{\Om}   | \p_t^j \Div_x \eta |^2.  \\
     \end{split}
\end{equation}
Likewise,
\begin{equation}\label{e-e-t-16}
     \begin{split}
         \int_{\Om} & \  \overline{b}^2  \Big\{  \p_3 (\p_t^j \eta)\cdot \p_3(\p_t^j \vu)-\Div_x(\p_t^j \eta) \p_3 (\p_t^j u_3)    \\
          & \   +  \Div_x(\p_t^j \eta) \Div_x(\p_t^j \vu)   -\p_3 (\p_t^j \eta_3) \Div_x(\p_t^j \vu) \Big\}       \\
        =  & \   \overline{b}^2 \f{1}{2} \f{\rm{d}}{\dt}\int_{\Om}  \left(
          | \p_t^j \p_3 \eta  |^2 -2 \Div_x (\p_t^j \eta) \p_3 (\p_t^j \eta_3) +| \Div_x (\p_t^j \eta) |^2
     \right)                  \\
        =  & \  \overline{b}^2 \f{1}{2} \f{\rm{d}}{\dt}\int_{\Om}  \left( | \p_t^j \p_3 \eta_{\ast}  |^2 +|\p_t^j \Div_{\ast}\eta_{\ast} |^2             \right).\\
     \end{split}
\end{equation}
To proceed, we estimate the right-hand side of (\ref{e-e-t-12}). By Lemma \ref{le-1},
\begin{equation}\label{e-e-t-17}
     \begin{split}
        \left| \int_{\Om} \left( J \p_t^j \vu \cdot \mathbf{N}^j+ \overline{\vt}^{-1} J \p_t^j \vt M^j\right)\right| & \ \lesssim
\|  \p_t^j \vu \|_0 \|\mathbf{N}^j \|_0+ \|  \p_t^j \vt \|_0 \|M^j \|_0 \\
    & \        \lesssim   \sqrt{\mathcal{D}_n} \sqrt{ \mathcal{E}_n  \mathcal{D}_n   }     \\
         & \ =    \sqrt{ \mathcal{E}_n    } \mathcal{D}_n.    \\
     \end{split}
\end{equation}
Other integrals on the right-hand side of (\ref{e-e-t-12}) are bounded through
\begin{equation}\label{e-e-t-18}
\|\p_t^j \eta  \|_1  \| \eta \|_3 \|\p_t^j \vu  \|_1  \lesssim
\sqrt{\mathcal{D}_n}    \sqrt{\mathcal{E}_n}  \sqrt{\mathcal{D}_n} = \sqrt{ \mathcal{E}_n    } \mathcal{D}_n.
\end{equation}

Based on (\ref{e-e-t-15})-(\ref{e-e-t-18}), we obtain from (\ref{e-e-t-12}), by summing over $j$ from $0$ to $n$ and recalling the definition of $ \mathcal{E}_n^{(t)}$, that
\begin{equation}\label{e-e-t-19}
     \begin{split}
     \f{1}{2} \f{\rm{d}}{\dt} \mathcal{E}_n^{(t)}   & \
+\sum_{j=0}^n \int_{\Om}  \Big(
\mu J |\Grad_{\ma}\p_t^j \vu|^2 \\
     & \  +(\mu+\lambda)J|\Div_{\ma}\p_t^j \vu|^2+\kappa \overline{\vt}^{-1} J |\Grad_{\ma}\p_t^j \vt|^2
\Big)   \\
      \lesssim    & \   \sqrt{ \mathcal{E}_n    } \mathcal{D}_n. \\
     \end{split}
\end{equation}
Observing that
\begin{equation}\label{e-e-t-20}
     \begin{split}
        J |\Grad_{\ma}\p_t^j \vu|^2= & \  |\Grad_x (\p_t^j \vu)|^2+ (J-1)  |\Grad_x (\p_t^j \vu)|^2\\
     & \      + J \left(    \Grad_{\ma}\p_t^j \vu+ \Grad_x (\p_t^j \vu)        \right)  :\left( \Grad_{\ma}\p_t^j \vu  -\Grad_x (\p_t^j \vu)       \right );              \\
     \end{split}
\end{equation}
\begin{equation}\label{e-e-t-21}
     \begin{split}
        J |\Grad_{\ma}\p_t^j \vt|^2= & \  |\Grad_x (\p_t^j \vt)|^2+ (J-1)  |\Grad_x (\p_t^j \vt)|^2\\
     & \      + J \left(    \Grad_{\ma}\p_t^j \vt+ \Grad_x (\p_t^j \vt)        \right)  :\left( \Grad_{\ma}\p_t^j \vt  -\Grad_x (\p_t^j \vt)       \right );              \\
     \end{split}
\end{equation}
\begin{equation}\label{e-e-t-22}
     \begin{split}
        J |\Div_{\ma}\p_t^j \vu|^2= & \  |\Div_x (\p_t^j \vu)|^2+ (J-1)  |\Div_x (\p_t^j \vu)|^2\\
     & \      + J \left(    \Div_{\ma}\p_t^j \vu+ \Div_x (\p_t^j \vu)        \right)  :\left( \Div_{\ma}\p_t^j \vu  -\Div_x (\p_t^j \vu)       \right ).             \\
     \end{split}
\end{equation}
We may invoke (\ref{e-e-t-20})-(\ref{e-e-t-22}) to deduce
\begin{equation}\label{e-e-t-23}
     \begin{split}
        \sum_{j=0}^n & \ \int_{\Om}  \Big(
\mu J |\Grad_{\ma}\p_t^j \vu|^2 +(\mu+\lambda)J|\Div_{\ma}\p_t^j \vu|^2 +\kappa \overline{\vt}^{-1} J |\Grad_{\ma}\p_t^j \vt|^2
\Big)  \\
         & \  \geq  \sum_{j=0}^n    \int_{\Om}\Big( \mu  |\Grad_x (\p_t^j \vu)|^2 +(\mu+\lambda) |\Div_x (\p_t^j \vu)|^2+ \kappa \overline{\vt}^{-1} |\Grad_x (\p_t^j \vt)|^2\Big)
         -C  \sqrt{ \mathcal{E}_n    } \mathcal{D}_n,  \\
     \end{split}
\end{equation}
where we have estimated the remaining terms from (\ref{e-e-t-20})-(\ref{e-e-t-22}) in the same manner as (\ref{e-e-t-18}) and Lemma \ref{le-1}. Finally, inserting (\ref{e-e-t-23}) into (\ref{e-e-t-19}), employing the
Poincar\'{e} inequality and recalling the definition of $\mathcal{D}_n^{(t)}$, we arrive at (\ref{e-e-t-11}). This finishes the proof of Lemma \ref{le-2}.            $\Box$

\subsection{Energy evolution of horizontal derivatives}\label{e-e-h}
The main purpose of this subsection is to derive the horizontal spatial derivatives of solutions. To this end, we adopt the following linear perturbed form by reformulating (\ref{e-e-t-2}) as
\begin{equation}\label{e-e-h-1}
\left\{\begin{aligned}
& \p_t \eta=\vu, \\
& \overline{\vr} \p_t \vu -\mu \Delta_{x}\vu -(\mu+\lambda)\Grad_{x}\Div_{x}\vu +  \overline{\vr} \Grad_{x} \vt  -  \overline{\vr}  \overline{\vt} \Grad_x \Div_x \eta   \\
&= \overline{b}^2 (   \p_{33}\eta -\p_3 \Div_x\eta \mathbf{e}_3 +\Grad_x \Div_x \eta- \Grad_x \p_3 \eta_3     ) +\mathbf{W},\\
& \overline{\vr} \p_t \vt+   \overline{\vr}\overline{\vt} \Div_{x} \vu
=\kappa \Delta_{x} \vt+V ,\\
& (\eta,\vu,\vt)|_{t=0}=(\eta_0,\vu_0,\vt_0),\\
& \vu|_{\p \Om}=\mathbf{0},\,\,\,  \vt |_{\p \Om}=0. \\
\end{aligned}\right.
\end{equation}
Here the nonlinear term $\mathbf{W}=(W_1,W_2,W_3)$ with
\begin{equation}\label{e-e-h-2}
     \begin{split}
         W_i= & \   R^{\eta}_i -   \overline{\vr}  (J^{-1}-1)\p_t u_i -  \overline{\vr}   (       \ma_{ik} -\delta_{ik}      )  \p_k \vt        \\
    & \  -  \overline{\vr}  \ma_{ik} \p_k (  (J^{-1}-1)\vt     )               \\
         & \  +\mu \ma_{jk} \p_k \ma_{jl}\p_l u_i  + \mu (   \ma_{jk} \ma_{jl}-\delta_{jk}\delta_{jl}    )  \p_k\p_l u_i            \\
         & \  +(\mu+\lambda) \ma_{ij} \p_j \ma_{kl} \p_l u_k+ (\mu+\lambda) (\ma_{ij} \ma_{kl} -\delta_{ij} \delta_{kl}    )  \p_j  \p_l u_k,\,\,\,\, i=1,2,3,
     \end{split}
\end{equation}
and
\begin{equation}\label{e-e-h-3}
     \begin{split}
         V= & \  -   \overline{\vr}  (J^{-1}-1)\p_t \vt-   \overline{\vr}   \overline{\vt}     (   \ma_{ik}\p_k u_i -\delta_{ik}\p_k u_i   )             \\
         & \  +\kappa \ma_{jk} \p_k \ma_{jl}\p_l \vt  + \kappa (   \ma_{jk} \ma_{jl}-\delta_{jk}\delta_{jl}    )  \p_k\p_l \vt    \\
         & \ -  \Big(
            \overline{\vr} \vt+ \overline{\vr} (J^{-1}-1)\vt + \overline{\vr} (J^{-1}-1) \overline{\vt}
            \Big) \ma_{ik}\p_k u_i \\
         & \  +\Big\{ \f{\mu}{2}  |     \ma_{ik}\p_k u_l +\ma_{lk}\p_k u_i               |^2 +\lambda (\ma_{ik}\p_k u_i)^2 \Big\}  .
     \end{split}
\end{equation}

For the nonlinear terms $\mathbf{W}$ and $V$, we estimate as follows. The proof is similar to Lemma \ref{le-1} and thus omitted, cf. also Lemma 3.3 in \cite{TW}.
\begin{Lemma}\label{le-3}
Let $ n \geq 3$ be an integer. It holds
\begin{equation}\label{e-e-h-4}
\| \overline{\Grad}^{2n-2} \mathbf{W} \|_0^2  +\| \overline{\Grad}^{2n-2}  V \|_0^2   \lesssim  (\mathcal{E}_n)^2,
\end{equation}
\begin{equation}\label{e-e-h-5}
\| \overline{\Grad}^{2n-1} \mathbf{W} \|_0^2  +\| \overline{\Grad}^{2n-1} V \|_0^2   \lesssim  \mathcal{E}_n  \mathcal{D}_n.
\end{equation}

\end{Lemma}

We introduce the horizontal energy as
\begin{equation}\label{e-e-h-6}
     \begin{split}
       \mathcal{E}_n^{(\ast)}= & \      \overline{\vr} \|\Grad_{\ast} \vu \|_{0,2n-1}^2 +\overline{\vr} \overline{\vt}^{-1} \|\Grad_{\ast} \vt \|_{0,2n-1}^2
       +  \overline{\vr}  \overline{\vt} \| \Grad_{\ast} \Div_x \eta   \|_{0,2n-1}^2   \\
     & \   +\overline{b}^2 \| \Grad_{\ast} \p_3 \eta_{\ast} \|_{0,2n-1}^2 +\overline{b}^2 \| \Grad_{\ast} \Div_{\ast} \eta_{\ast} \|_{0,2n-1}^2 ,  \\
     \end{split}
\end{equation}
and the horizontal dissipation as
\begin{equation}\label{e-e-h-7}
 \mathcal{D}_n^{(\ast)}=
\| \Grad_{\ast} \vu \|_{1,2n-1}^2+\| \Grad_{\ast} \vt \|_{1,2n-1}^2.
\end{equation}
Now we give the energy evolution of horizontal spatial derivatives. In analogy with Lemma \ref{le-2}, the cancellations between the momentum equation and temperature equation is useful.
\begin{Lemma}\label{le-4}
Let $ n \geq 3$ be an integer. Then
\begin{equation}\label{e-e-h-8}
\f{\rm{d}}{\dt} \mathcal{E}_n^{(\ast)} +\mathcal{D}_n^{(\ast)} \lesssim
\sqrt{\mathcal{E}_n}  \mathcal{D}_n.
\end{equation}

\end{Lemma}
{\bf{Proof.} } For any given $\alpha \in \mathbb{N}^2$ with $1\leq |\alpha|\leq 2n$. We first apply the differential operator $\p^{\alpha}$ to (\ref{e-e-h-1})$_2$ and (\ref{e-e-h-1})$_3$ respectively\footnote{In fact, we should write the horizontal differential operator as $\p_{1}^{\alpha_1}\p_{2}^{\alpha_2}$ with $ \alpha=(\alpha_1,\alpha_2)$. Nevertheless, here and in what follows we choose the simplified notation $\p^{\alpha}$ only for brevity.   }, then take $L^2$ inner products with $\p ^{\alpha}\vu$ and $\overline{\vt}^{-1}\p ^{\alpha}\vt$ respectively, finally integrate by parts to deduce that
\begin{equation}\label{e-e-h-9}
     \begin{split}
          \f{1}{2} \f{\rm{d}}{\dt} & \   \int_{\Om} \Big(    \overline{\vr} |\p ^{\alpha}\vu|^2 + \overline{\vr}\overline{\vt}^{-1} |\p ^{\alpha}\vt|^2  \\
        & \  +  \overline{\vr}\overline{\vt} |\p^{\alpha} \Div_x \eta|^2 +\overline{b}^2 |\p^{\alpha} \p_3 \eta_{\ast}|^2+   \overline{b}^2 |\p^{\alpha}  \Div_{\ast}\eta_{\ast}|^2 \Big)                      \\
         & \ +   \int_{\Om} \Big( \mu | \Grad_x \p^{\alpha}\vu |^2 +(\mu+\lambda) |\Div_x \p^{\alpha}\vu|^2+ \kappa \overline{\vt}^{-1} | \Grad_x \p^{\alpha}\vt |^2           \Big)               \\
        = &  \int_{\Om}\Big( \p^{\alpha} \mathbf{W} \cdot \p^{\alpha}\vu +\overline{\vt}^{-1}  \p^{\alpha} V \cdot \p^{\alpha} \vt \Big).
     \end{split}
\end{equation}
To proceed, we estimate the right-hand side of (\ref{e-e-h-9}). In accordance with the assumption of $\alpha$, we make the decomposition $\alpha=(\alpha-\beta)+\beta$ with $\beta\in \mathbb{N}^2,|\beta|=1$. Thanks to
(\ref{e-e-h-5}) and integration by parts,
\begin{equation}\label{e-e-h-10}
     \begin{split}
        \int_{\Om} & \   \Big( \p^{\alpha} \mathbf{W} \cdot \p^{\alpha}\vu +\overline{\vt}^{-1}  \p^{\alpha} V \cdot \p^{\alpha} \vt \Big)\\
      & \    =-  \int_{\Om}  \Big( \p^{\alpha-\beta}\mathbf{W} \cdot \p^{\alpha+\beta}\vu + \overline{\vt}^{-1} \p^{\alpha-\beta}V \cdot \p^{\alpha+\beta}  \vt\Big)               \\
         & \   \lesssim   \| \p^{\alpha-\beta}\mathbf{W} \|_0  \| \p^{\alpha+\beta}\vu \|_0 +\|  \p^{\alpha-\beta}V \|_0    \| \p^{\alpha+\beta}  \vt \|_0    \\
         & \  \lesssim  \| \mathbf{W} \|_{2n-1} \| \vu\|_{2n+1}+ \| V \|_{2n-1} \| \vt\|_{2n+1} \\
         & \   \lesssim \sqrt{ \mathcal{E}_n\mathcal{D}_n}\sqrt{\mathcal{D}_n}= \sqrt{ \mathcal{E}_n}\mathcal{D}_n.
     \end{split}
\end{equation}
We thus obtain (\ref{e-e-h-8}) upon combining (\ref{e-e-h-9})-(\ref{e-e-h-10}), applying the Poincar\'{e} inequality and running over all $\alpha\in \mathbb{N}^2$ such that$1\leq |\alpha|\leq 2n$.   $\Box$

\subsection{Energy evolution of the flow mapping}\label{e-e-f}
In this subsection, we derive the estimates on dissipations of the flow mapping $\eta$. It should be stressed that the Lagrangian formulation of the MHD system is crucial. To begin with, we report the following technical lemma from \cite{TW}, cf. Lemma 3.5 therein.
\begin{Lemma}\label{le-5}
Let $N\geq 4$ be an integer. Then
\begin{equation}\label{e-e-f-1}
\| \mathbf{W} \|^2_{4N-1} \lesssim  (\mathcal{D}_{2N})^2  + \mathcal{E}_{N+2}\mathcal{E}_{2N},
\end{equation}
\begin{equation}\label{e-e-f-2}
\| \mathbf{W} \|^2_{2(N+2)} \lesssim   \mathcal{E}_{2N}\mathcal{D}_{N+2}.
\end{equation}

\end{Lemma}

Following \cite{TW}, we define the recovering energy as
\begin{equation}\label{e-e-f-3}
 \mathcal{E}_n^{ ( \sharp) } = \mu \| \Grad_x \eta \|^2_{0,2n}  +(\mu+\lambda) \| \Div_x \eta \| ^2_{0,2n}
\end{equation}
and the recovering dissipation as
\begin{equation}\label{e-e-f-4}
 \mathcal{D}_n^{  (\sharp) } = \| \Div_x \eta \| ^2_{0,2n}+   \| \p_3 \eta \| ^2_{0,2n} + \|\eta\|^2_{0,2n}.
\end{equation}

The energy evolution of the flow mapping read as
\begin{Lemma}\label{le-6}
Let $N\geq 4$ be an integer. Then
\begin{equation}\label{e-e-f-5}
     \begin{split}
         \f{\rm{d}}{\dt} & \    \left\{  \mathcal{E}_{2N}^{  (\sharp) } +2 \sum_{\alpha\in \mathbb{N}^2,|\alpha|\leq 4N}\int_{\Om}  \overline{\vr}\p^{\alpha}\vu  \cdot \p^{\alpha}\eta   \right\} +  \mathcal{D}_{2N}^{  (\sharp) } \\
       & \  \lesssim  \sqrt{   \mathcal{E}_{2N}} \mathcal{D}_{2N} +\sqrt{   \mathcal{E}_{N+2}} \mathcal{E}_{2N} +   \mathcal{D}_{2N}^{ (t) }  +   \mathcal{D}_{2N}^{ (\ast) }  ,      \\
     \end{split}
\end{equation}
\begin{equation}\label{e-e-f-6}
     \begin{split}
         \f{\rm{d}}{\dt} & \    \left\{  \mathcal{E}_{N+2}^{  (\sharp) } +2 \sum_{\alpha\in \mathbb{N}^2,|\alpha|\leq 2(N+2)}\int_{\Om}  \overline{\vr}\p^{\alpha}\vu  \cdot \p^{\alpha}\eta  \right\} +  \mathcal{D}_{N+2}^{  (\sharp) }   \\
       & \  \lesssim  \sqrt{   \mathcal{E}_{2N}} \mathcal{D}_{N+2}  +   \mathcal{D}_{N+2}^{ (t) }  +   \mathcal{D}_{N+2}^{ (\ast) }  .      \\
     \end{split}
\end{equation}

\end{Lemma}
{\bf{Proof.} } For brevity, we denote by $n$ for either $2N$ or $N+2$. Let $\alpha \in \mathbb{N}^2$ with $|\alpha| \leq 2n$. We first apply the differential operator $\p^{\alpha}$ to (\ref{e-e-h-1})$_2$ and then take $L^2$ inner product with $\p ^{\alpha}\eta$, finally integrate by parts to obtain
\begin{equation}\label{e-e-f-7}
     \begin{split}
       \int_{\Om}  \overline{\vr}\p_t(\p^{\alpha}\vu)  \cdot \p^{\alpha}\eta   & \ + \f{1}{2}\f{\rm{d}}{\dt}  \int_{\Om} \Big(  \mu  | \Grad_x \p^{\alpha}\eta|^2  +(\mu+\lambda)|\Div_x \p^{\alpha}\eta |^2   \Big)             \\
       & \  +\int_{\Om} \Big( \overline{\vr}\overline{\vt} |\p^{\alpha} \Div_x \eta|^2 +\overline{b}^2 |\p^{\alpha} \p_3 \eta_{\ast}  |^2  +\overline{b}^2 |\p^{\alpha} \Div_{\ast}\eta_{\ast}    |^2  \Big)           \\
         =& \ \int_{\Om} \p^{\alpha}\eta \cdot \p^{\alpha}\mathbf{W}   - \int_{\Om} \overline{\vr}  \Grad_x (\p^{\alpha}\vt)\cdot \p^{\alpha}\eta.            \\
     \end{split}
\end{equation}
With the help of (\ref{e-e-h-1})$_1$, we simplify the first integral on the left-hand side of (\ref{e-e-f-7}) as
\begin{equation}\label{e-e-f-8}
 \int_{\Om}  \overline{\vr}\p_t(\p^{\alpha}\vu)  \cdot \p^{\alpha}\eta   =  \f{\rm{d}}{\dt}   \int_{\Om}  \overline{\vr}\p^{\alpha}\vu  \cdot \p^{\alpha}\eta
 -      \int_{\Om}  \overline{\vr} |\p^{\alpha}\vu|^2 .
\end{equation}
Obviously, it holds
\begin{equation}\label{e-e-f-9}
\int_{\Om}  \overline{\vr} |\p^{\alpha}\vu|^2
\lesssim
 \mathcal{D}_{n}^{ (t) }  +   \mathcal{D}_{n}^{ (\ast) }.
\end{equation}
Next, we estimate the second integral on the right-hand side of (\ref{e-e-f-7}). If $|\alpha|=0$,
\begin{equation}\label{e-e-f-10}
     \begin{split}
         \left| -\overline{\vr}\int_{\Om}  \Grad_x \vt \cdot \eta  \right| & \ \lesssim \| \Grad_x \vt \|_0  \| \eta \|_0    \\
     & \  \lesssim   \sqrt{  \mathcal{D}_{n}^{ (t) }   }     \sqrt{  \mathcal{D}_{n}^{ (\sharp) }   }          \\
         & \  \lesssim  \ep \mathcal{D}_{n}^{ (\sharp) }+ \f{C}{\ep}  \mathcal{D}_{n}^{ (t) } ;            \\
     \end{split}
\end{equation}
If $0<|\alpha|\leq 2n $,
\begin{equation}\label{e-e-f-11}
     \begin{split}
         \left| \int_{\Om} \overline{\vr}  \Grad_x (\p^{\alpha}\vt)\cdot \p^{\alpha}\eta  \right| & \ \lesssim \| \Grad_x (\p^{\alpha}\vt) \|_0  \| \p^{\alpha}\eta \|_0    \\
     & \  \lesssim   \| \p^{\alpha}\vt \|_1  \| \p^{\alpha}\eta\|_0         \\
         & \  \lesssim   \sqrt{  \mathcal{D}_{n}^{ (\ast) }   }     \sqrt{  \mathcal{D}_{n}^{ (\sharp) }   }              \\
          & \  \lesssim  \ep \mathcal{D}_{n}^{ (\sharp) }+ \f{C}{\ep}  \mathcal{D}_{n}^{ (\ast) } .           \\
     \end{split}
\end{equation}
Here $\ep>0$ is a parameter to be determined later. Thus we are left to estimate the first integral on the right-hand side of (\ref{e-e-f-7}). If $n=N+2$,
\begin{equation}\label{e-e-f-12}
     \begin{split}
         \left| \int_{\Om} \p^{\alpha}\eta \cdot \p^{\alpha}\mathbf{W}  \right| & \ \lesssim \| \p^{\alpha}\eta \|_0  \| \p^{\alpha}\mathbf{W} \|_0    \\
     & \  \lesssim   \sqrt{  \mathcal{D}_{N+2}  }     \sqrt{  \mathcal{E}_{2N}  \mathcal{D}_{N+2}  } ,         \\
     \end{split}
\end{equation}
where (\ref{e-e-f-2}) was used. If $n=2N$ and $|\alpha| \leq 4N-1$,
\begin{equation}\label{e-e-f-13}
     \begin{split}
         \left| \int_{\Om} \p^{\alpha}\eta \cdot \p^{\alpha}\mathbf{W}  \right| & \ \lesssim \| \p^{\alpha}\eta \|_0  \| \mathbf{W} \|_{4N-1}    \\
     & \  \lesssim   \sqrt{  \mathcal{D}_{2N}  }     \sqrt{  \mathcal{E}_{2N}  \mathcal{D}_{2N}  } ,         \\
     \end{split}
\end{equation}
where (\ref{e-e-h-5}) was invoked. If $n=2N$ and $|\alpha| =4N$, we make the decomposition $\alpha=(\alpha-\beta)+\beta$ for some $\beta \in \mathbb{N}^2,|\beta|=1$. In view of (\ref{e-e-f-1}) and integration by parts,
\begin{equation}\label{e-e-f-14}
     \begin{split}
         \int_{\Om} \p^{\alpha}\eta \cdot \p^{\alpha}\mathbf{W}  & \  =-  \int_{\Om} \p^{\alpha+\beta}\eta \cdot \p^{\alpha-\beta} \mathbf{W}      \\
     & \  \lesssim   \|  \eta \|_{4N+1}\|  \mathbf{W} \|_{4N-1} ,         \\
     & \  \lesssim   \sqrt{  \mathcal{E}_{2N}  } (        \sqrt{  \mathcal{D}_{2N}  } +\sqrt{  \mathcal{D}_{N+2}    \mathcal{E}_{2N}  }       ). \\
     \end{split}
\end{equation}

Therefore, combining (\ref{e-e-f-8})-(\ref{e-e-f-14}) and fixing $\ep>0$ sufficiently small, we conclude from (\ref{e-e-f-7}) the relations (\ref{e-e-f-5})-(\ref{e-e-f-6}) upon summing over all $\alpha \in \mathbb{N}^2,|\alpha|\leq 2n$ and applying the Poincar\'{e} inequality.          $\Box$

\section{Improved estimates}\label{i-e}
We are now in a position to derive the vertical derivatives of solutions and improve the energy evolution estimates with the help of elliptic estimates.

\subsection{Energy evolution of vertical derivatives}\label{e-e-v}
As observed by Matsumura and Nishida \cite{MaNi1} when treating the full Navier-Stokes system in domains with boundary, it is indispensable to separate the momentum equation as the vertical part and the horizontal part if one expects to derive the vertical derivative estimates. For brevity, we introduce the quantity $q:=-\Div_x \eta$ which satisfies $\p_t q=-\Div_x \vu$. We may take the third component of the momentum equation (\ref{e-e-h-1})$_2$ to find that
\begin{equation}\label{e-e-v-1}
(2\mu+\lambda)\p_t\p_3 q+ \overline{\vr}\overline{\vt} \p_3 q= - \overline{\vr} \p_t u_3 + \mu \Delta_{\ast} u_3-\mu \p_3 \Div_{\ast} \vu_{\ast} -\overline{\vr} \p_3 \vt+ W_3,
\end{equation}
and similarly take the horizontal components as
\begin{equation}\label{e-e-v-2}
     \begin{split}
         -\mu \p_{33} \vu_{\ast} - \overline{b}^2  \p_{33}\eta_{\ast} = & \       - \overline{\vr} \p_t \vu_{\ast}  +\mu \Delta_{\ast} \vu_{\ast} -(\mu+\lambda) \Grad_{\ast} \p_t q-    \overline{\vr}\overline{\vt}    \Grad_{\ast} q                    \\
         & \   -  \overline{b}^2  ( \Grad_{\ast} q +\Grad_{\ast} \p_3 \eta_3) -  \overline{\vr} \Grad_{\ast} \vt  +\mathbf{W}_{\ast}.            \\
     \end{split}
\end{equation}
It turns out that (\ref{e-e-v-1}) and (\ref{e-e-v-2}) obey the ODE structure. More precisely, they could be written in a unified way as
\[
c_i \f{\rm{d}}{\dt}  f_i +c_i f_i =g_i,  \,\,\, i=1,2,
\]
with $f_i$ being $\p_3 q$ and $ \p_{33}\eta_{\ast}$ respectively; $c_i$ are constants and $g_i$ are the source terms.

Below we give the energy evolution of vertical derivatives for the flow mapping. Although the proof follows the same line as \cite{TW}, we shall present it for convenience.
\begin{Lemma}\label{le-7}
Let $n\geq 3$ be an integer. It holds that
\begin{equation}\label{e-e-v-3}
     \begin{split}
        \f{\rm{d}}{\dt}  \mathcal{V}_n  +& \   \| \vu \|^2_{2n+1} +\| \Div_x \eta \|^2_{2n}   + \| \p_3 \eta \|^2_{2n}  +\|  \eta \| ^2_{2n}                \\
        & \             \lesssim   \mathcal{E}_n  \mathcal{D}_n +   \mathcal{D}_{n}^{ (t) }  +   \mathcal{D}_{n}^{ (\ast) }  +   \mathcal{D}_{n}^{ (\sharp) } +\|\p_t \vu  \|^2_{2n-1}  +\|  \vt \|_{2n}^2 ,                \\
     \end{split}
\end{equation}
where $\mathcal{V}_n $ is equivalent to $\| \p_3 q \|_{2n-1}^2+\| \p_{33}\eta_{\ast} \|_{2n-1}^2$.
\end{Lemma}
{\bf{Proof.} } Let $0\leq k\leq 2n-1$ be fixed. By taking the norm $\| \cdot \|_{k,2n-k-1}^2$ on both sides of (\ref{e-e-v-1}),
\begin{equation}\label{e-e-v-4}
     \begin{split}
         (2\mu+\lambda)\overline{\vr}\overline{\vt}\f{\rm{d}}{\dt} \| \p_3 q \|^2_{k,2n-k-1}       & \  + \overline{\vr}^2\overline{\vt}^2 \| \p_3 q  \| ^2_{k,2n-k-1}
         + (2\mu+\lambda)^2  \| \p_3 \p_t q \|    ^2_{k,2n-k-1}                         \\
          & \  =\|    - \overline{\vr} \p_t u_3 + \mu \Delta_{\ast} u_3-\mu \p_3 \Div_{\ast} \vu_{\ast} -\overline{\vr} \p_3 \vt+ W_3   \|   ^2_{k,2n-k-1}                                 \\
         & \ \lesssim   \|\p_t u_3   \|^2_{2n-1} + \| \vu \| ^2_{k+1,2n-k}   + \| \vt\|_{2n}^2    + \|  W_3   \| ^2_{2n-1} .                 \\
     \end{split}
\end{equation}

Observing that
\begin{equation}\label{e-e-v-5}
     \begin{split}
         \| \p_{33} \eta_3\|^2_{k,2n-k-1}  & \  = \| \p_3(q+ \Div_{\ast}\eta_{\ast}) \|^2_{k,2n-k-1}               \\
          & \  \lesssim   \| q  \|^2_{k+1,2n-k-1} +\| \p_3 \eta_{\ast} \| ^2_{k,2n-k};             \\
     \end{split}
\end{equation}
\begin{equation}\label{e-e-v-6}
     \begin{split}
         \| \p_{33} u_3\|^2_{k,2n-k-1}  & \  = \| \p_3(\p_t q+ \Div_{\ast}\vu_{\ast}) \|^2_{k,2n-k-1}               \\
          & \  \lesssim   \| \p_t q  \|^2_{k+1,2n-k-1} +\| \vu_{\ast} \| ^2_{k+1,2n-k}.              \\
     \end{split}
\end{equation}
Based on (\ref{e-e-v-4})-(\ref{e-e-v-6}),
\begin{equation}\label{e-e-v-7}
     \begin{split}
        \f{\rm{d}}{\dt}\| \p_3 q \|^2_{k,2n-k-1}    & \    +\| q  \|^2_{k+1,2n-k-1} +\| \p_t q  \|^2_{k+1,2n-k-1}           \\
                 & \  + \| \p_{33} u_3\|^2_{k,2n-k-1} +  \| \p_{33} \eta_3\|^2_{k,2n-k-1}        \\
          \lesssim  & \ \|\p_t u_3   \|^2_{2n-1}+ \| \vu \| ^2_{k+1,2n-k} +\| \p_3 \eta_{\ast} \| ^2_{k,2n-k} +  \| \vt\|_{2n}^2    + \|  W_3   \| ^2_{2n-1}             \\
         & \ +  \mathcal{D}_{n}^{ (t) }  +   \mathcal{D}_{n}^{ (\ast) }  +   \mathcal{D}_{n}^{ (\sharp) } .    \\
     \end{split}
\end{equation}

In a similar manner, we take the norm $\| \cdot \|_{k,2n-k-1}^2$ on both sides of (\ref{e-e-v-2}) to conclude that
\begin{equation}\label{e-e-v-8}
     \begin{split}
         \mu \overline{b}^2 \f{\rm{d}}{\dt} \|  \p_{33}\eta_{\ast}\|^2_{k,2n-k-1}    & \ + \overline{b}^4 \|\p_{33}\eta_{\ast}  \| ^2_{k,2n-k-1}
         +\mu^2   \|\p_{33}\vu_{\ast}  \| ^2_{k,2n-k-1}                         \\
       =& \  \|   - \overline{\vr} \p_t \vu_{\ast}  +\mu \Delta_{\ast} \vu_{\ast} -(\mu+\lambda) \Grad_{\ast} \p_t q-    \overline{\vr}\overline{\vt}    \Grad_{\ast} q           \\
         & \  -  \overline{b}^2  ( \Grad_{\ast} q +\Grad_{\ast} \p_3 \eta_3) -  \overline{\vr} \Grad_{\ast} \vt  +\mathbf{W}_{\ast}  \|   ^2_{k,2n-k-1}       \\
          \lesssim & \ \|\p_t \vu_{\ast}  \|^2_{2n-1}+ \|  \vu_{\ast} \|^2_{k,2n-k+1}+ \|\p_t q  \|^2_{k,2n-k} + \| q  \|^2_{k,2n-k}  \\
          & \  + \|  \p_3 \eta_3 \|^2_{k,2n-k} +\| \vt \|^2_{2n}+ \| \mathbf{W}_{\ast}  \|^2_{k,2n-1}. \\
     \end{split}
\end{equation}
It follows from (\ref{e-e-v-7})-(\ref{e-e-v-8}) that
\begin{equation}\label{e-e-v-9}
     \begin{split}
             \f{\rm{d}}{\dt}      & \   \left(  \| \p_3 q \|^2_{k,2n-k-1}  +  \|  \p_{33}\eta_{\ast}\|^2_{k,2n-k-1}        \right)        \\
         & \  +      \| \p_t q  \|^2_{k+1,2n-k-1}+\| q  \|^2_{k+1,2n-k-1}  +\| \p_3 \eta  \|^2_{k+1,2n-k-1} + \| \vu  \|^2_{k+2,2n-k-1}                \\
       \lesssim  & \   \| \p_t q  \|^2_{k,2n-k}+\| q  \|^2_{k,2n-k}  +\| \p_3 \eta  \|^2_{k,2n-k} + \| \vu  \|^2_{k+1,2n-k}         \\
         & \  + \| \p_t \vu \|^2_{2n-1} +\|  \vt\|^2_{2n}+ \| \mathbf{W}  \|^2_{2n-1} \\
         & \ +  \mathcal{D}_{n}^{ (t) }  +   \mathcal{D}_{n}^{ (\ast) }  +   \mathcal{D}_{n}^{ (\sharp) } .             \\
     \end{split}
\end{equation}
From the recursive inequality (\ref{e-e-v-9}) we deduce that there exist suitable constants $a_k>0,k=0,...,2n-1$ such that
\begin{equation}\label{e-e-v-10}
     \begin{split}
         \f{\rm{d}}{\dt}   \sum_{k=0}^{2n-1}      & \  a_k   \left(  \| \p_3 q \|^2_{k,2n-k-1}  +  \|  \p_{33}\eta_{\ast}\|^2_{k,2n-k-1}        \right)     \\
     & \  +   \sum_{k=0}^{2n-1}   \Big(  \| \p_t q  \|^2_{k+1,2n-k-1}+\| q  \|^2_{k+1,2n-k-1}                  \\
     & \   +\| \p_3 \eta  \|^2_{k+1,2n-k-1} + \| \vu  \|^2_{k+2,2n-k-1}   \Big)  \\
         \lesssim  & \   \|\p_t q \|^2_{0,2n}  +\| q \|^2_{0,2n}   +\|\vu \|^2_{1,2n}  +\|\p_3 \eta \|^2_{0,2n}               \\
         & \  + \| \p_t \vu \|^2_{2n-1} +\|  \vt\|^2_{2n}+ \| \mathbf{W}  \|^2_{2n-1}  \\
         & \ +  \mathcal{D}_{n}^{ (t) }  +   \mathcal{D}_{n}^{ (\ast) }  +   \mathcal{D}_{n}^{ (\sharp) }             \\
       \lesssim  & \    \| \p_t \vu \|^2_{2n-1} +\|  \vt\|^2_{2n}+ \| \mathbf{W}  \|^2_{2n-1}  \\
        & \ +  \mathcal{D}_{n}^{ (t) }  +   \mathcal{D}_{n}^{ (\ast) }  +   \mathcal{D}_{n}^{ (\sharp) } .            \\
     \end{split}
\end{equation}
Finally, we conclude from (\ref{e-e-v-10}) that
\begin{equation}\label{e-e-v-11}
     \begin{split}
        \f{\rm{d}}{\dt} \mathcal{V}_n   & \  + \|  \p_t q \|^2_{2n} + \|  q \|^2_{2n} + \|  \p_3 \eta \|^2_{2n}+ \|  \vu \|^2_{2n+1}     \\
      \lesssim  & \    \| \p_t \vu \|^2_{2n-1} +\|  \vt\|^2_{2n}+ \| \mathbf{W}  \|^2_{2n-1}  \\
        & \ +  \mathcal{D}_{n}^{ (t) }  +   \mathcal{D}_{n}^{ (\ast) }  +   \mathcal{D}_{n}^{ (\sharp) } ,           \\
     \end{split}
\end{equation}
where
\[
\mathcal{V}_n:= \sum_{k=0}^{2n-1}a_k  \left(  \| \p_3 q \|^2_{k,2n-k-1}  +  \|  \p_{33}\eta_{\ast}\|^2_{k,2n-k-1}        \right).
\]
Therefore, one obtains (\ref{e-e-v-3}) from (\ref{e-e-v-11}) immediately upon noticing $\| \mathbf{W}  \|^2_{2n-1}\leq \mathcal{E}_n\mathcal{D}_n$ by (\ref{e-e-h-5}).       $\Box$

\subsection{Energy improvements by elliptic estimates}\label{e-i-e}
In this subsection, we improve the energy estimates of $(\vu,\vt)$ in terms of elliptic regularity theory. First of all, we rewrite the momentum equation (\ref{e-e-h-1})$_2$ as
\begin{equation}\label{e-i-e-1}
\left\{\begin{aligned}
&  -\mu \Delta_{x}\vu -(\mu+\lambda)\Grad_{x}\Div_{x}\vu   = -\overline{\vr} \p_t \vu  -  \overline{\vr} \Grad_{x} \vt  +  \overline{\vr}  \overline{\vt} \Grad_x \Div_x \eta   \\
&+\overline{b}^2 (   \p_{33}\eta -\p_3 \Div_x\eta \mathbf{e}_3 +\Grad_x \Div_x \eta- \Grad_x \p_3 \eta_3     ) +\mathbf{W},\\
& \vu|_{\p \Om}=\mathbf{0}, \\
\end{aligned}\right.
\end{equation}
and the temperature equation (\ref{e-e-h-1})$_3$ as
\begin{equation}\label{e-i-e-2}
\left\{\begin{aligned}
&  -\kappa \Delta_{x} \vt =-        \overline{\vr} \p_t \vt-   \overline{\vr}\overline{\vt} \Div_{x} \vu
+V ,\\
& \vt |_{\p \Om}=0. \\
\end{aligned}\right.
\end{equation}
We then apply the classical elliptic regularity theory to prove
\begin{Lemma}\label{le-8}
Let $n \geq 3$ be an integer. Then
\begin{equation}\label{e-i-e-3}
     \begin{split}
         \sum_{j=1}^n \| \p_t^j \vu \|^2_{2n-2j+1}     & \ + \sum_{j=0}^n \| \p_t^j \vt \|^2_{2n-2j+1}  \\
      & \  \lesssim    \mathcal{E}_n \mathcal{D}_n +    \mathcal{D}_{n}^{ (t) }           +\|  \vu \|^2_{2n}.             \\
     \end{split}
\end{equation}
\end{Lemma}
{\bf{Proof.} } For any given $j\in \{ 1,...,n-1 \}$, we apply the time derivative $\p_t^j $ to (\ref{e-i-e-1})$_1$ and invoke the well-known elliptic regularity theory (see for instance \cite{ADN}) to conclude that
\begin{equation}\label{e-i-e-4}
     \begin{split}
         \| \p_t^j \vu \|^2_{2n-2j+1}   \lesssim   & \    \| \p_t^{j+1} \vu \|^2_{2n-2j-1}+ \|  \Grad_x^2  \p_t^{j} \eta \|^2_{2n-2j-1}        \\
         & \  +\| \Grad_x \p_t^j \vt \|  ^2_{2n-2j-1}  +
            \|  \p_t^j \mathbf{W} \|  ^2_{2n-2j-1}     \\
         \lesssim   & \      \| \p_t^{j+1} \vu \|^2_{2n-2(j+1)-1}+ \|   \p_t^{j} \eta \|^2_{2n-2j+1}                 \\
         & \  + \|  \p_t^j \vt \|  ^2_{2n-2j} +\| \overline{\Grad}^{2n-1} \mathbf{W} \|_0^2.\\
     \end{split}
\end{equation}
It follows that
\begin{equation}\label{e-i-e-5}
     \begin{split}
         \sum_{j=1}^n \| \p_t^j \vu \|^2_{2n-2j+1}    \lesssim    & \  \| \p_t^n \vu\|^2_{1} +  \sum_{j=1}^{n-1} \| \p_t^j \eta \|^2_{2n-2j+1}                            \\
         & \   +  \sum_{j=1}^{n-1} \|  \p_t^j \vt \|  ^2_{2n-2j} +  \| \overline{\Grad}^{2n-1} \mathbf{W} \|_0^2       \\
        \lesssim    & \ \sum_{j=1}^{n-1} \|  \p_t^{j-1} \vu \|  ^2_{2n-2(j-1)-1}    +  \sum_{j=1}^{n-1} \|  \p_t^j \vt \|  ^2_{2n-2j}           \\
         & \ +  \mathcal{D}_{n}^{ (t) }  + \| \overline{\Grad}^{2n-1} \mathbf{W} \|_0^2 \\
         =  & \      \sum_{j=0}^{n-2} \|  \p_t^{j} \vu \|  ^2_{2n-2j-1}    +  \sum_{j=1}^{n-1} \|  \p_t^j \vt \|  ^2_{2n-2j}                     \\
         & \ +  \mathcal{D}_{n}^{ (t) }  + \| \overline{\Grad}^{2n-1} \mathbf{W} \|_0^2; \\
     \end{split}
\end{equation}
whence with the help of Sobolev interpolation inequality and Young inequality we furthermore obtain
\begin{equation}\label{e-i-e-6}
     \begin{split}
         \sum_{j=1}^n \| \p_t^j \vu \|^2_{2n-2j+1}    \lesssim    & \      \| \vu \|^2_{2n-1} +   \sum_{j=1}^{n-2} \|  \p_t^{j} \vu \|  ^2_{2n-2j-1}    +  \sum_{j=1}^{n-1} \|  \p_t^j \vt \|  ^2_{2n-2j}                           \\
         & \ +  \mathcal{D}_{n}^{ (t) }  + \| \overline{\Grad}^{2n-1} \mathbf{W} \|_0^2 \\
        \lesssim   & \   \| \vu \|^2_{2n-1} + \sum_{j=1}^{n-2} \|  \p_t^{j} \vu \|  ^2_{0}   +  \sum_{j=1}^{n-1} \|  \p_t^j \vt \|  ^2_{2n-2j}  \\
        & \   +  \mathcal{D}_{n}^{ (t) }  + \| \overline{\Grad}^{2n-1} \mathbf{W} \|_0^2 \\
           \lesssim  & \    \| \vu \|^2_{2n-1}+  \sum_{j=1}^{n-1} \|  \p_t^j \vt \|  ^2_{2n-2j}+  \mathcal{D}_{n}^{ (t) }  + \| \overline{\Grad}^{2n-1} \mathbf{W} \|_0^2.  \\
     \end{split}
\end{equation}

In a similar manner, for any given $j\in \{ 0,1,...,n-1 \}$, we apply the time derivative $\p_t^j $ to (\ref{e-i-e-2})$_1$ and invoke the well-known elliptic regularity theory (see for instance \cite{ADN}) to conclude that
\begin{equation}\label{e-i-e-7}
     \begin{split}
        \| \p_t^j \vt \|^2_{2n-2j+1} \lesssim   & \  \| \p_t^{j+1} \vt \|^2_{2n-2j-1} +\| \Div_x \p_t ^j \vu  \|^2_{2n-2j-1} +\| \p_t^j V \|^2_{2n-2j-1}      \\
    \lesssim   & \  \| \p_t^{j+1} \vt \|^2_{2n-2(j+1)-1}  +  \| \p_t ^j \vu  \|^2_{2n-2j}  +\| \p_t^j V \|^2_{2n-2j-1}.       \\
     \end{split}
\end{equation}
The recursive inequality (\ref{e-i-e-7}) particularly gives rise to
\begin{equation}\label{e-i-e-8}
     \begin{split}
       \sum_{j=0}^n   \| \p_t^j \vt \|^2_{2n-2j+1} \lesssim   & \  \| \p_t^n  \vt \|^2_{1} +  \sum_{j=0}^{n-1}    \| \p_t ^j \vu  \|^2_{2n-2j} +\| \overline{\Grad}^{2n-1} V \|_0^2     \\
      =  & \   \| \p_t^n  \vt \|^2_{1}  + \|\vu  \|^2_{2n} +  \sum_{j=1}^{n-1}    \| \p_t ^j \vu  \|^2_{2n-2j} +\| \overline{\Grad}^{2n-1} V \|_0^2.  \\
     \end{split}
\end{equation}
Combining (\ref{e-i-e-6}) and (\ref{e-i-e-8}), we make use of Sobolev interpolation inequality and Young inequality again to deduce that
\begin{equation}\label{e-i-e-9}
     \begin{split}
       \sum_{j=1}^n \| \p_t^j \vu \|^2_{2n-2j+1} +\sum_{j=0}^n   \| \p_t^j \vt \|^2_{2n-2j+1}\lesssim   & \    \| \vu \|^2_{2n-1}  +\sum_{j=1}^{n-1} \|  \p_t^j \vt \|  ^2_{0}       \\
    & \ +  \| \p_t^n  \vt \|^2_{1}  + \|\vu  \|^2_{2n}   +  \sum_{j=1}^{n-1}    \| \p_t ^j \vu  \|^2_{0}      \\
         & \  +  \mathcal{D}_{n}^{ (t) }  + \| \overline{\Grad}^{2n-1} \mathbf{W} \|_0^2 +\| \overline{\Grad}^{2n-1} V \|_0^2  \\
      \lesssim   & \   \mathcal{E}_n \mathcal{D}_n +    \mathcal{D}_{n}^{ (t) }           +\|  \vu \|^2_{2n}, \\
     \end{split}
\end{equation}
where we used (\ref{e-e-h-5}) and the definition of $ \mathcal{D}_{n}^{ (t) }$ in the last step. This completes the proof of Lemma \ref{le-8}.      $\Box$

In the next lemma, we control the energy $\mathcal{E}_n$ suitably, which plays a crucial role in deriving the global estimates.
\begin{Lemma}\label{le-9}
Let $n \geq 3$ be an integer. Then
\begin{equation}\label{e-i-e-10}
\mathcal{E}_n   \lesssim
\mathcal{E}_{n}^{ (t) } +  \mathcal{E}_{n}^{ (\sharp) }    +\mathcal{V}_n +(\mathcal{E}_n )^2 .
\end{equation}
\end{Lemma}
{\bf{Proof.} } For any given $j\in \{0, 1,...,n-1 \}$, we apply the time derivative $\p_t^j $ to (\ref{e-i-e-1})$_1$ and invoke the well-known elliptic regularity theory (see for instance \cite{ADN}) to conclude that
\begin{equation}\label{e-i-e-11}
     \begin{split}
         \| \p_t^j \vu \|^2_{2n-2j}   \lesssim   & \    \| \p_t^{j+1} \vu \|^2_{2n-2j-2}+ \|  \Grad_x^2  \p_t^{j} \eta \|^2_{2n-2j-2}        \\
         & \  +\| \Grad_x \p_t^j \vt \|  ^2_{2n-2j-2}  +
            \|  \p_t^j \mathbf{W} \|  ^2_{2n-2j-2}     \\
         \lesssim   & \      \| \p_t^{j+1} \vu \|^2_{2n-2(j+1)}+ \|   \p_t^{j} \eta \|^2_{2n-2j}                 \\
         & \  + \|  \p_t^j \vt \|  ^2_{2n-2j-1} +\| \overline{\Grad}^{2n-2} \mathbf{W} \|_0^2.\\
     \end{split}
\end{equation}
We deduce from the recursive inequality (\ref{e-i-e-11}) that
\begin{equation}\label{e-i-e-12}
     \begin{split}
         \sum_{j=0}^n  \| \p_t^j \vu \|^2_{2n-2j}  \lesssim  & \  \| \p_t^n \vu \|_0^2 +  \sum_{j=0}^{n-1}  \| \p_t^j \eta \|^2_{2n-2j}                       \\
     & \  +  \sum_{j=0}^{n-1} \|  \p_t^j \vt \|  ^2_{2n-2j-1}  +\| \overline{\Grad}^{2n-2} \mathbf{W} \|_0^2        \\
         \lesssim & \  \mathcal{E}_{n}^{ (t) } +  \|\eta \|^2_{2n}  +    \sum_{j=1}^{n-1}   \| \p_t^{j-1} \vu  \|^2_{2n-2j}            \\
         & \   +  \sum_{j=0}^{n-1} \|  \p_t^j \vt \|  ^2_{2n-2j-1}  +\| \overline{\Grad}^{2n-2} \mathbf{W} \|_0^2        \\
           \lesssim  & \    \sum_{j=0}^{n-2}   \| \p_t^{j} \vu  \|^2_{2n-2j-2}   +\sum_{j=0}^{n-1} \|  \p_t^j \vt \|  ^2_{2n-2j-1}  \\
           & \   + \mathcal{E}_{n}^{ (t) } + \mathcal{E}_{n}^{ (\sharp) } +\mathcal{V}_n+ \| \overline{\Grad}^{2n-2} \mathbf{W} \|_0^2, \\
     \end{split}
\end{equation}
where we have used the estimate $\|\eta\|^2_{2n+1}  \lesssim \mathcal{E}_{n}^{ (\sharp) } +  \mathcal{V}_n $ in the last step. Analogously, for any given $j\in \{ 0,1,...,n-1 \}$, we apply the time derivative $\p_t^j $ to (\ref{e-i-e-2})$_1$ and invoke the well-known elliptic regularity theory (see for instance \cite{ADN}) to conclude that
\begin{equation}\label{e-i-e-13}
     \begin{split}
        \| \p_t^j \vt \|^2_{2n-2j} \lesssim   & \  \| \p_t^{j+1} \vt \|^2_{2n-2j-2} +\| \Div_x \p_t ^j \vu  \|^2_{2n-2j-2} +\| \p_t^j V \|^2_{2n-2j-2}      \\
    \lesssim   & \  \| \p_t^{j+1} \vt \|^2_{2n-2(j+1)}  +  \| \p_t ^j \vu  \|^2_{2n-2j-1}  +\| \overline{\Grad}^{2n-2} V \|_0^2,      \\
     \end{split}
\end{equation}
yielding recursively
\begin{equation}\label{e-i-e-14}
     \begin{split}
       \sum_{j=0}^n  \| \p_t^j \vt \|^2_{2n-2j} \lesssim   & \  \| \p_t^{n} \vt \|^2_{0} +    \sum_{j=0}^{n-1}     \|  \p_t ^j \vu  \|^2_{2n-2j-1} +\| \overline{\Grad}^{2n-2} V \|_0^2      \\
    \lesssim   & \ \mathcal{E}_{n}^{ (t) }    +    \sum_{j=0}^{n-1}     \|  \p_t ^j \vu  \|^2_{2n-2j-1} +\| \overline{\Grad}^{2n-2} V \|_0^2.      \\
     \end{split}
\end{equation}
Taking (\ref{e-i-e-12}) and (\ref{e-i-e-14}) into account, we infer by Sobolev interpolation inequality and Young inequality that
\begin{equation}\label{e-i-e-15}
     \begin{split}
         \sum_{j=0}^n  \| \p_t^j \vu \|^2_{2n-2j} +   \sum_{j=0}^n  \| \p_t^j \vt \|^2_{2n-2j} \lesssim  & \  \sum_{j=0}^{n-1}   \| \p_t^{j} \vu  \|^2_{2n-2j-2}   +\sum_{j=0}^{n-1} \|  \p_t^j \vt \|  ^2_{2n-2j-1} \\
     & \   + \mathcal{E}_{n}^{ (t) } + \mathcal{E}_{n}^{ (\sharp) } +\mathcal{V}_n  \\
         & \  +     \| \overline{\Grad}^{2n-2} \mathbf{W} \|_0^2+\| \overline{\Grad}^{2n-2} V \|_0^2     \\
     \lesssim     & \       \sum_{j=0}^{n-1}   \| \p_t^{j} \vu  \|^2_{0}   +\sum_{j=0}^{n-1} \|  \p_t^j \vt \|  ^2_{0}               \\
     & \            +           \mathcal{E}_{n}^{ (t) } + \mathcal{E}_{n}^{ (\sharp) } +\mathcal{V}_n               \\
     & \  +     \| \overline{\Grad}^{2n-2} \mathbf{W} \|_0^2+\| \overline{\Grad}^{2n-2} V \|_0^2     \\
     \lesssim & \   \mathcal{E}_{n}^{ (t) } + \mathcal{E}_{n}^{ (\sharp) } +\mathcal{V}_n + (\mathcal{E}_n)^2, \\
     \end{split}
\end{equation}
where we utilized (\ref{e-e-h-4}) and the definition of $\mathcal{E}_{n}^{ (t) }$ in the last step. This finishes the proof of Lemma \ref{le-9} upon recalling the estimate $\|\eta\|^2_{2n+1}  \lesssim \mathcal{E}_{n}^{ (\sharp) } +  \mathcal{V}_n $.                      $\Box$

\section{Global energy estimates}\label{g-e}
In this section we shall glue all the energy evolution estimates and the improved estimates obtained in Sections \ref{e-e}-\ref{i-e} so as to prove (\ref{ma_th-4}).

\begin{Lemma}\label{le-10}
Let $n$ be either $2N$ or $N+2$. There exists a modified energy $\mathcal{E}^{mod}_n$ which is equivalent to $\mathcal{E}_n$ and satisfy
\begin{equation}\label{g-e-1}
\f{\rm{d}}{\dt}
\mathcal{E}^{mod}_{2N} + \mathcal{D}_{2N}  \lesssim
\sqrt{ \mathcal{E}_{N+2} }\mathcal{E}_{2N},
\end{equation}
\begin{equation}\label{g-e-2}
\f{\rm{d}}{\dt}
\mathcal{E}^{mod}_{N+2} + \mathcal{D}_{N+2}  \lesssim
0.
\end{equation}
\end{Lemma}
{\bf{Proof.} } It follows from (\ref{e-e-v-3}) and (\ref{e-i-e-3}) that
\begin{equation}\label{g-e-3}
\f{\rm{d}}{\dt}         \mathcal{V}_n +\mathcal{D}_n   \lesssim
\mathcal{E}_n \mathcal{D}_n +\| \vu \|^2_{2n}+\| \vt \|^2_{2n} +    \mathcal{D}_{n}^{ (t) }  +    \mathcal{D}_{n}^{ (\ast) }  +    \mathcal{D}_{n}^{ (\sharp) }  ;
\end{equation}
whence we apply Sobolev interpolation inequality and Young inequality to infer that
\begin{equation}\label{g-e-4}
     \begin{split}
        \f{\rm{d}}{\dt}         \mathcal{V}_n +\mathcal{D}_n   \lesssim  & \  \mathcal{E}_n \mathcal{D}_n +\| \vu \|^2_{0}+\| \vt \|^2_{0} +    \mathcal{D}_{n}^{ (t) }  +    \mathcal{D}_{n}^{ (\ast) }  +    \mathcal{D}_{n}^{ (\sharp) } \\
   \lesssim  & \   \mathcal{E}_n \mathcal{D}_n +  \mathcal{D}_{n}^{ (t) }  +    \mathcal{D}_{n}^{ (\ast) }  +    \mathcal{D}_{n}^{ (\sharp) }.              \\
     \end{split}
\end{equation}
For notational convenience, we define $\mathcal{Y}_n$ with $\mathcal{Y}_{2N}=\sqrt{   \mathcal{E}_{2N}} \mathcal{D}_{2N} +\sqrt{   \mathcal{E}_{N+2}} \mathcal{E}_{2N}$ and $\mathcal{Y}_{N+2}=\sqrt{   \mathcal{E}_{2N}} \mathcal{D}_{N+2}$. Combining Lemmas \ref{le-2}, \ref{le-4}, \ref{le-6} and (\ref{g-e-4}),
\begin{equation}\label{g-e-5}
     \begin{split}
         \f{\rm{d}}{\dt}    & \  \left\{  \mathcal{E}_{n}^{ (t) }  +    \mathcal{E}_{n}^{ (\ast) }            +\ep \left(    \mathcal{E}_{2N}^{  (\sharp) } +2 \sum_{\alpha\in \mathbb{N}^2,|\alpha|\leq 4N}\int_{\Om}  \overline{\vr}\p^{\alpha}\vu  \cdot \p^{\alpha}\eta      \right)        +\ep^2 \mathcal{V}_n   \right\} \\
      & \  +      \mathcal{D}_{n}^{ (t) }  +    \mathcal{D}_{n}^{ (\ast) }  +   \ep \mathcal{D}_{n}^{ (\sharp) }    +\ep^2 \mathcal{D}_n        \\
       \lesssim   & \     \mathcal{Y}_n  +\ep ( \mathcal{D}_{n}^{ (t) }  +    \mathcal{D}_{n}^{ (\ast) })   +\ep^2 (     \mathcal{D}_{n}^{ (t) }  +    \mathcal{D}_{n}^{ (\ast) }+ \mathcal{D}_{n}^{ (\sharp) }  )  ,                \\
     \end{split}
\end{equation}
with $\ep>0$ being a positive parameter. Hence we choose $\ep$ to be suitably small to conclude that
\begin{equation}\label{g-e-6}
     \begin{split}
         \f{\rm{d}}{\dt}    & \  \left\{  \mathcal{E}_{n}^{ (t) }  +    \mathcal{E}_{n}^{ (\ast) }            +\ep \left(    \mathcal{E}_{2N}^{  (\sharp) } +2 \sum_{\alpha\in \mathbb{N}^2,|\alpha|\leq 4N}\int_{\Om}  \overline{\vr}\p^{\alpha}\vu  \cdot \p^{\alpha}\eta      \right)        +\ep^2 \mathcal{V}_n   \right\} \\
      & \  +      \mathcal{D}_{n}^{ (t) }  +    \mathcal{D}_{n}^{ (\ast) }  +   \ep \mathcal{D}_{n}^{ (\sharp) }    +\ep^2 \mathcal{D}_n        \\
       \lesssim   & \     \mathcal{Y}_n   .           \\
     \end{split}
\end{equation}
Upon setting
\begin{equation}\label{g-e-7}
\mathcal{E}^{mod}_n:= \mathcal{E}_{n}^{ (t) }  +    \mathcal{E}_{n}^{ (\ast) }            +\ep \left(    \mathcal{E}_{2N}^{  (\sharp) } +2 \sum_{\alpha\in \mathbb{N}^2,|\alpha|\leq 4N}\int_{\Om}  \overline{\vr}\p^{\alpha}\vu  \cdot \p^{\alpha}\eta      \right)        +\ep^2 \mathcal{V}_n,
\end{equation}
one obtains from Lemma \ref{le-9} that
\begin{equation}\label{g-e-8}
\mathcal{E}_{n}  \lesssim
   +  \mathcal{E}^{mod}_n+ (\mathcal{E}_{n})^2;
\end{equation}
whence we find that $\mathcal{E}^{mod}_n$ is equivalent to $\mathcal{E}_{n}$ due to our a priori assumption $\mathcal{E}_{2N}(T)\leq \delta$ for sufficiently small $\delta>0$. Obviously (\ref{g-e-6}) yields (\ref{g-e-1})-(\ref{g-e-2}). $\Box$

As stressed in the Introduction, the key issues in the applications of two-tier energy method is to show first the boundedness of energy in high regularity norms $\mathcal{E}_{2N}$ and then the decay of energy in lower ones $\mathcal{E}_{N+2}$.

\begin{Lemma}\label{le-11}
Let $N \geq 4$ be an integer. Then
\begin{equation}\label{g-e-9}
\mathcal{E}_{2N}(t) +\int_0^t \mathcal{D}_{2N}(s) \ds \lesssim  \mathcal{E}_{2N}(0),\,\,\,\, t\in [0,T].
\end{equation}
\end{Lemma}
{\bf{Proof.} } An integration of (\ref{g-e-1}) with respect to time gives rise to
\begin{equation}\label{g-e-10}
     \begin{split}
       \mathcal{E}_{2N}(t)+ \int_0^t \mathcal{D}_{2N}(\tau)  \rm{d}\tau     & \  \lesssim       \mathcal{E}_{2N}(0) +\int_0^t  \sqrt{ \mathcal{E}_{N+2} }\mathcal{E}_{2N}(\tau)   \rm{d}\tau        \\
     & \  \lesssim    \mathcal{E}_{2N}(0) + \sup_{0\leq \tau \leq t}  \mathcal{E}_{2N}(\tau)   \int_0^t  \sqrt{\delta} (1+\tau)^{-N+2} \ds         \\
         & \ \lesssim    \mathcal{E}_{2N}(0) +  \sqrt{\delta} \sup_{0\leq \tau \leq t}  \mathcal{E}_{2N}(\tau),     \\
     \end{split}
\end{equation}
where in the second step we used the decay estimate $\mathcal{E}_{N+2} (\tau)  \leq  \delta (1+\tau)^{-2N+4}  $ thanks to our assumption $\mathcal{G}_{2N}(T)\leq \delta$. (\ref{g-e-9}) follows from (\ref{g-e-10}) readily upon fixing $\delta>0$ sufficiently small.        $\Box$

In the end, we prove the decay estimates in lower regularity norms; this is indispensable in closing the global estimates.
\begin{Lemma}\label{le-12}
Let $N \geq 4$ be an integer. Then
\begin{equation}\label{g-e-11}
\mathcal{E}_{N+2}(t) \lesssim  \mathcal{E}_{2N}(0) (1+t)^{-2N+4},\,\,\,\, t\in [0,T].
\end{equation}
\end{Lemma}
{\bf{Proof.} } As observed in \cite{GuTi1,GuWa1,SrGu1}, the central issue is the application of Sobolev interpolation inequality. More precisely, it holds
\begin{equation}\label{g-e-12}
     \begin{split}
        \| \Grad_{\ast}^{   2(N+2)+1   } \eta \|_0^2  \lesssim   & \     \| \Grad_{\ast}^{   2(N+2)  } \eta \|_0^{2 \beta}         \| \Grad_{\ast}^{   4N+1  } \eta \|_0^{2 (1-\beta) }                               \\
   \lesssim   & \   ( \mathcal{D}_{N+2}   )^{\beta}   ( \mathcal{E}_{2N}   )^{1-\beta}   ,           \\
     \end{split}
\end{equation}
with $\beta:= \f{2N-4}{2N-3}$. Recalling the definition of $\mathcal{E}_n,\mathcal{D}_n$, we conclude
\begin{equation}\label{g-e-13}
\mathcal{E}_{N+2}
\lesssim
( \mathcal{D}_{N+2}   )^{\beta}   ( \mathcal{E}_{2N}   )^{1-\beta},
\end{equation}
which together with (\ref{g-e-9}) implies
\begin{equation}\label{g-e-14}
\mathcal{E}^{mod}_{N+2}   \lesssim  \mathcal{E}_{N+2} \lesssim  [\mathcal{E}_{2N} (0)]^{1-\beta}( \mathcal{D}_{N+2}   )^{\beta}.
\end{equation}
We then substitute the relation (\ref{g-e-14}) into (\ref{g-e-2}) so as to obtain
\begin{equation}\label{g-e-15}
 \f{\rm{d}}{\dt}    \mathcal{E}^{mod}_{N+2} +  \f{C}{    [\mathcal{E}_{2N} (0)]^{\gamma}     } (\mathcal{E}^{mod}_{N+2})^{1+\gamma} \leq 0,
\end{equation}
with $\gamma:= \f{1}{2N-4}  $. Finally, we solve the above differential inequality to conclude
\begin{equation}\label{g-e-16}
\mathcal{E}_{N+2}(t)  \lesssim  \mathcal{E}^{mod}_{N+2}  (t) \lesssim  \f{   \mathcal{E}_{2N} (0)     }
{      \Big(       [\mathcal{E}_{2N} (0)]^{\gamma}     +  \gamma C [\mathcal{E}_{N+2} (0)]^{\gamma}        t   \Big) ^{\f{1}{\gamma}}               }\mathcal{E}_{N+2} (0),\,\,\,\,t\in [0,T].
\end{equation}
It follows that
\begin{equation}\label{g-e-17}
\mathcal{E}_{N+2}(t)  \lesssim   \f{  \mathcal{E}_{2N} (0)    }{   1+t^{  \f{1}{\gamma}         }           } =\f{  \mathcal{E}_{2N} (0)    }{   1+t^{  2N-4        }           } .
\end{equation}
This completes the proof of Lemma \ref{le-12}.       $\Box$

Consequently, we get the global energy estimate (\ref{ma_th-4}) by Lemmas \ref{le-11}, \ref{le-12}. The proof of our main Theorem \ref{m-th} is thus completely finished.

\section{Concluding remarks}

\begin{itemize}
\item{
We consider the three-dimensional \emph{full} compressible viscous non-resistive MHD system. Global well-posedness is proved in an infinite slab $\R^2 \times (0,1)$ with strong background magnetic field. In order to capture the weak dissipations of magnetic field in the non-resistive fluids, we appeal to the Lagrangian coordinates, inspired by the incompressible viscous non-resistive MHD system \cite{LXZ,XZ} and the barotropic compressible MHD system \cite{TW}. The main tool is the celebrated two-tier energy method developed by Guo and Tice in a series of work of viscous surface waves without surface tension \cite{GuTi1,GuTi2}. To our best knowledge, this seems to be the first \emph{global-in-time} result for full compressible viscous non-resistive MHD system in multi-dimensions.
}

\item{
It should be noticed that our main result relies heavily on the geometric domain and the background magnetic field. More precisely, the finite deepness of domain is useful when applying the Poincar\'{e} inequality, cf.
Lemmas \ref{le-2}, \ref{le-4}. For the background magnetic field $\overline{\vb}=(\overline{B}_1,\overline{B}_2,\overline{B}_3)^{t}$, we require that $\overline{B}_3 \neq 0$. This reveals in some sense that the vertical magnetic field enjoys the stabilization effects in full compressible viscous non-resistive MHD. Indeed, our previous work \cite{LS2D,LS2Dfu} on global weak solutions of two-dimensional compressible non-resistive MHD system confirms also this crucial fact; here the velocity field takes the form $(u_1,u_2,0)(t,x_1,x_2)$, while the magnetic field takes the form $(0,0,b)(t,x_1,x_2)$.
}

\end{itemize}

{\large{{\bf{\centerline{Acknowledgements}}}}}
The research of Y. Li is supported by National Natural Science Foundation of China under grant number 12001003. The author thanks Prof. Y. Sun for constant guidance and encouragement.

%\abcdefg

\end{document}